\documentclass[invmat,final,numbook,francais]{svjour}
\usepackage{times}
\usepackage[T1]{fontenc}
\usepackage{amsmath}
\usepackage{amssymb}
\usepackage{amsfonts}
\usepackage{eucal}
\usepackage[all]{xy}
\usepackage[frenchb]{babel}
\usepackage{pslatex}

\newtheorem{prop}{Proposition}[subsection]
\newtheorem{theo}[prop]{Théor\`eme}
\newtheorem{coro}[prop]{Corollaire}
\newtheorem{lemm}[prop]{Lemme}
\newtheorem{lemm*}[equation]{Lemme}

\spnewtheorem{conj}[prop]{Conjecture}{\bf}{\rm}
\spnewtheorem{vide}[prop]{}{\bf}{\rm}
\spnewtheorem{defi}[prop]{Définition}{\bf}{\rm}

\spnewtheorem{rema}[prop]{Remarques}{\bf}{\rm}
\spnewtheorem{exem}[prop]{Exemple}{\bf}{\rm}
\spnewtheorem{nota}[prop]{Notations}{\bf}{\rm}

\numberwithin{equation}{prop}

\newcommand{\riso}{ \tilde {\rightarrow}\, }
\newcommand{\liso}{ \tilde {\leftarrow}\, }

\renewcommand{\sp}{\mathrm{sp}}

\newcommand{\FF}{{\mathcal{F}}}

\newcommand{\E}{{\mathcal{E}}}

\newcommand{\D}{{\mathcal{D}}}

\newcommand{\PP}{{\mathcal{P}}}

\renewcommand{\O}{{\mathcal{O}}}

\newcommand{\V}{\mathcal{V}}
\renewcommand{\S}{\mathcal{S}}
\newcommand{\Y}{\mathcal{Y}}

\newcommand{\X}{\mathfrak{X}}

\newcommand{\U}{\mathfrak{U}}

\newcommand{\A}{\mathbb{A}}
\renewcommand{\P}{\mathbb{P}}

\newcommand{\DD}{\mathbb{D}}
\renewcommand{\L}{\mathbb{L}}
\newcommand{\R}{\mathbb{R}}
\newcommand{\Q}{\mathbb{Q}}
\newcommand{\Z}{\mathbb{Z}}

\newcommand{\hdag}{  \phantom{}{^{\dag} }    }

\begin{document}

\title{$F$-isocristaux surconvergents et surcohérence différentielle}
\titlerunning{$F$-isocristaux surconvergents et surcohérence différentielle}

\author{D. Caro
\thanks{L'auteur a bénéficié du soutien du réseau européen TMR \textit{Arithmetic Algebraic Geometry}
(contrat numéro UE MRTN-CT-2003-504917).}
}                     
%
%
\institute{Université Paris-Sud, Laboratoire de Mathématiques d'Orsay,
Orsay Cedex 91405, France, \email{daniel.caro@math.u-psud.fr}}
%

%
\maketitle
\begin{abstract}
Soient $\V$ un anneau de valuation discrète complet d'inégales caractéristiques,
de corps résiduel parfait $k$,
$\PP$ un $\V$-schéma formel propre et lisse, $T$ un diviseur de la fibre spéciale $P$ de $\PP$,
  $U$ l'ouvert de $P$ complémentaire de $T$,
  $Y$ un sous-$k$-schéma fermé lisse de $U$.
Nous prouvons que
  la catégorie des $F$-isocristaux surconvergents sur $Y$ est équivalente à celle des
  $F$-isocristaux surcohérents sur $Y$ (voir \cite[6.2.1 et 6.4.3.a)]{caro_devissge_surcoh}).
  Plus généralement, nous établissons par recollement une telle équivalence pour tout $k$-schéma séparé lisse $Y$.
  Nous vérifions de plus que les $F$-complexes de $\D ^\dag _{\PP} (\hdag T) _\Q$-modules à cohomologie bornée
  et $\D ^\dag _{\PP} (\hdag T) _\Q$-surcohérente se dévissent en $F$-isocristaux surconvergents.

\end{abstract}


\tableofcontents


\section*{Introduction}
Soit $k$ un corps parfait de caractéristique $p>0$.
Nous savons depuis les années 60 que la constructibilité en cohomologie étale $l$-adique sur les
$k$-variétés algébriques (i.e. $k$-schémas séparés de type fini) est, lorsque $l \not =p$,
stable par les {\og six opérations cohomologiques de Grothendieck\fg},
à savoir $\otimes$,  $\mathcal{H}om$, $f _*$, $f ^*$, $f _!$ et $f^!$.
Nous souhaitons un analogue $p$-adique, i.e.,
disposer d'une cohomologie $p$-adique stable et définie sur les $k$-variétés algébriques.
Or, les complexes de $\D$-modules (par défaut à gauche)
à cohomologie bornée, holonome et régulière sur les
variétés algébriques sur {\it un corps de caractéristique nulle} sont stables par ces six opérations.
Berthelot et, par une approche voisine
(on dispose de théorèmes de comparaison entre les deux constructions, voir \cite{huyghe-comparaison}),
Mebkhout et Narvaez-Macarro (voir \cite{Mebkhout-Narvaez-Macarro_D})
ont alors construit une version $p$-adique de la théorie des $\D$-modules, i.e., {\og les $\D$-modules arithmétiques\fg}
(voir \cite{Be1}, \cite{Be2}, \cite{Beintro2}).
Nous utiliserons exclusivement dans ces travaux les $\D$-modules arithmétiques au sens de Berthelot, même si la version
de Mebkhout et Narvaez-Macarro a été une source d'inspiration dans notre travail.
À l'instar de la caractéristique nulle,
certains $F$-complexes (le $F$ indique que l'objet est muni d'une structure de Frobenius)
de $\D$-modules arithmétiques devraient constitués les coefficients $p$-adiques
d'une cohomologie $p$-adique stable (voir par exemple \cite{caro_surholonome}).
La construction des $\D$-modules arithmétiques s'est fortement inspirée de celle des isocristaux surconvergents,
coefficients $p$-adiques de la cohomologie rigide.
Cette cohomologie rigide représente
une unification des deux précédentes cohomologies $p$-adiques : celle cristalline et celle de Dwork-Monsky-Washnitzer
(voir \cite{Becris}, \cite{MonWas68}, \cite{Berig_ini} et \cite{Berig}).
Les $F$-isocristaux surconvergents
correspondent à un analogue $p$-adique des $\Q _l$-faisceaux constructibles lisses (on bénéficie ainsi de la stabilité par
image inverse mais nous ne disposons pas de notion d'images directes de $F$-isocristaux surconvergents),
tandis que les $F$-complexes de $\D$-modules arithmétiques contiennent un analogue $p$-adique des
$\Q_l$-faisceaux constructibles.

Après avoir replacé notre travail dans son contexte, expliquons en gros ses principaux résultats.
Nous avions vérifié dans \cite{caro_devissge_surcoh}
que pour toute $k$-variété lisse $Y$,
il existe un ouvert $\smash{\widetilde{Y}}$ affine et dense de $Y$
tel que la catégorie
des $F$-isocristaux surconvergents sur $\smash{\widetilde{Y}}$
soit incluse dans celle des $F$-$\D$-modules arithmétiques sur $\smash{\widetilde{Y}}$.
Nous prolongeons ici cette inclusion sur $Y$ tout entier.
La propriété essentielle caractérisant cette image est la {\og surcohérence différentielle\fg}
(pour la notion de surcohérence, on se référera à \cite[3]{caro_surcoherent}).
Réciproquement, nous vérifierons que les $F$-complexes de $\D $-modules arithmétiques
à cohomologie bornée et surcohérente se dévissent en $F$-isocristaux surconvergents.
Comme la surcohérence est stable par image directe par un morphisme propre, image inverse extraordinaire,
foncteur cohomologique local à support strict dans un sous-schéma fermé
(la stabilité par les autres opérations reste encore à établir),
nous avons ainsi englobé la catégorie
des $F$-isocristaux surconvergents sur $Y$ dans une catégorie notamment stable par image directe par un morphisme propre
(ou autrement dit, l'image directe par un morphisme propre d'un $F$-isocristal surconvergent sur $Y$
n'est pas à priori un $F$-isocristal surconvergent mais un
$F$-complexe de $\D $-modules arithmétiques
à cohomologie bornée et surcohérente).

Précisons à présent le contenu de ce papier.
Soient $\V$ un anneau de valuation discrète complet d'inégales caractéristiques,
de corps résiduel parfait $k$ de caractéristique $p>0$,
$\PP$ un $\V$-schéma formel propre et lisse, $T$ un diviseur de la fibre spéciale $P$ de $\PP$,
  $\U$ l'ouvert de $\PP$ complémentaire de $T$,
  $Y$ un sous-$k$-schéma fermé lisse de la fibre spéciale $U$ de $\U$,
  $\overline{Y}$ l'adhérence de $Y$ dans $P$.
  Nous avions défini dans \cite[6.2.1 et 6.4.3.a)]{caro_devissge_surcoh}
  une sous-catégorie pleine de celle des $F$-$\D ^\dag _{\PP} (\hdag T) _\Q$-modules surcohérents
   constituée par
  {\og les $F$-isocristaux surcohérents sur $Y$\fg}.
  Dans la première partie, lorsque $Y$ est en outre affine, nous prouvons que
  la catégorie des $F$-isocristaux surconvergents sur $Y$ est équivalente à celle des
  $F$-isocristaux surcohérents sur $Y$ (voir \ref{eqcataffine}).
Afin d'étendre cette équivalence de catégories au cas où $Y$ est une $k$-variété lisse quelconque
(qui ne se plonge plus forcément dans un $\V$-schéma formel propre et lisse etc.),
nous donnons, dans cette même partie, une procédure de recollement
que l'on utilisera à la fois pour construire la catégorie des $F$-isocristaux surconvergents sur $Y$ et
celle des $F$-isocristaux surcohérents sur $Y$ (voir \ref{isocsurcvgen-defi}).
La vérification de la procédure de recollement concernant les $F$-isocristaux surconvergents fait l'objet
de la première partie (voir \ref{ind-isocsurcvgen-vide}) tandis que
celle concernant les $F$-isocristaux surcohérents
est traitée dans la seconde partie (voir \ref{ind-isosurcoh}).
On obtient alors par recollement l'équivalence de catégories voulue dans le cas général (voir \ref{eqcat-gen}).
Cette équivalence entraîne que la catégorie
des $F\text{-}\D ^\dag _{\PP} (\hdag T) _\Q$-modules surcohérents $\E$ à support dans $\overline{Y}$ tels que
$\E | \U$ soit dans l'image essentielle de $\sp _{Y \hookrightarrow \U, +}$
(pour la description de cette image essentielle, voir \cite[4.1.9]{caro-construction})
est équivalente à celle des $F$-isocristaux surcohérents sur $Y$ (voir \ref{surcoh=dagdag}).
Dans la dernière partie, nous en déduirons que les $F$-complexes de $\D ^\dag _{\PP} (\hdag T) _\Q$-modules
à cohomologie bornée et $\D ^\dag _{\PP} (\hdag T) _\Q$-surcohérente sont dévissables en $F$-isocristaux surconvergents (voir \ref{caradev}),
ce qui généralise \cite[8.1.4-8.1.5]{caro_devissge_surcoh}.
En étendant la notion de dévissabilité en $F$-isocristaux surconvergents au cas des $F$-complexes quasi-cohérents
(nous conservons les notations de \cite[1.1.3]{caro_courbe-nouveau}),
nous construisons ensuite la sous-catégorie pleine
$F \text{-} \smash[b]{\underset{^{\longrightarrow }}{LD }}  ^\mathrm{b} _{\Q, \textrm{dév}}
(\overset{^\mathrm{g}}{} \smash{\widehat{\D}} _{\PP} ^{(\bullet)} (T ))$ de
$F \text{-} \smash[b]{\underset{^{\longrightarrow }}{LD }}  ^\mathrm{b} _{\Q, \mathrm{qc}}
(\overset{^\mathrm{g}}{} \smash{\widehat{\D}} _{\PP} ^{(\bullet)} (T ))$ des
$F$-complexes dévissables
en $F$-isocristaux surconvergents (voir \ref{def-dev-genU}).
Cette définition est compatible avec celle de la sous-catégorie pleine
$F \text{-}D ^\mathrm{b} _{\textrm{dév}}
(\smash{\D} ^\dag _{\PP} (\hdag T) _\Q)$
de
$F\text{-}D ^\mathrm{b} _{\textrm{coh}}
(\smash{\D} ^\dag _{\PP} (\hdag T) _\Q)$
des $F$-complexes dévissables en $F$-isocristaux surconvergents
que nous avions construite en
\cite[8.1.1]{caro_devissge_surcoh} (voir \ref{dev=qwcdev+coh}).
Nous vérifions enfin que
$F \text{-} \smash[b]{\underset{^{\longrightarrow }}{LD }}  ^\mathrm{b} _{\Q, \textrm{dév}}
(\overset{^\mathrm{g}}{} \smash{\widehat{\D}} _{\PP} ^{(\bullet)} (T ))$
est la plus petite sous-catégorie pleine triangulée de
$F \text{-} \smash[b]{\underset{^{\longrightarrow }}{LD }}  ^\mathrm{b} _{\Q, \textrm{qc}}
(\overset{^\mathrm{g}}{} \smash{\widehat{\D}} _{\PP} ^{(\bullet)} (T ))$
contenant les $F$-isocristaux surconvergents sur les sous-variétés lisses de $P$ (voir \ref{ppetitetrgl}).

\section*{Notations}
Soit $\V$ un anneau de valuation discrète complet,
de corps résiduel parfait $k$ de caractéristique $p>0$, de corps de
fractions $K$ de caractéristique $0$.
De plus, $s\geq 1$ sera un entier et $F$ la puissance
$s$-ème de l'endomorphisme de Frobenius. Les modules sont par défaut des modules à gauche.
Si $\E$ est un faisceau abélien, on pose $\E _\Q:=\E \otimes _\Z \Q$.

En général,
les $\V$-schémas formels faibles seront désignés par des lettres romanes surmontées du symbole {\og $\dag$\fg},
les $\V$-schémas formels par des lettres calligraphiques ou gothiques,
les $k$-schémas par des lettres romanes, e.g.,
$X ^\dag$, $\X$, $X$.
De plus, si $X ^\dag$ est un $\V$-schéma formel faible,
$X $ désignera sa fibre spéciale et $\X$ le $\V$-schéma formel égal à son complété $p$-adique etc.

Si $f$ : $\PP ' \rightarrow \PP$ est un morphisme de $\V$-schémas formels lisses,
 $f _0$ : $P' \rightarrow P$ sera le morphisme induit (et de même pour les morphismes de
 $\V$-schémas formels faibles lisses).
 Par abus de notations, on pourra encore l'indiquer par $f$.
Lorsque $T$ est un diviseur de $P$ tel que $f ^{-1} (T)$ soit un diviseur de $P'$
nous désignerons par
$f ^! _{T}$ ou $f ^! $
l'image inverse extraordinaire par $f$ à singularités surconvergentes le long de $T$
et par $f _{T,+}$ ou $f _+$ l'image directe par $f$ à singularités surconvergentes le long de $T$
(voir \cite[3.4, 3.5, 4.3]{Beintro2} et \cite[1.1.5]{caro_courbe-nouveau}).
De plus, si $Z$ est un sous-schéma fermé de $P$,
le foncteur cohomologique local
à support strict dans $Z$
et le foncteur restriction en dehors de $Z$
sont respectivement notés $\R \underline{\Gamma} ^\dag _Z $ et $(\hdag Z)$ (voir \cite[2.2.6]{caro_surcoherent}).
Ces foncteurs sont bien définis sur la catégorie des ($F$-)complexes quasi-cohérents
$(F\text{-})\smash{\underset{^{\longrightarrow}}{LD}} ^{\mathrm{b}} _{\Q ,\mathrm{qc}}
(\smash{\widehat{\D}} _{\PP} ^{(\bullet)}(T))$ (notations de \cite[1.1.3]{caro_courbe-nouveau}).
Si $\E \in (F\text{-})\smash{\underset{^{\longrightarrow}}{LD}} ^{\mathrm{b}} _{\Q ,\mathrm{qc}}
(\smash{\widehat{\D}} _{\PP} ^{(\bullet)}(T))$, on utilisera aussi la notation
$\E (\hdag Z)$ à la place de $(\hdag Z) (\E)$.
Pour tout diviseur $T$ de $P$, nous désignerons par $\DD _{\PP,T}$ ou $\DD _T$,
le dual $\smash{\D} ^\dag _{\PP} (\hdag T) _{\Q}$-linéaire (voir \cite[I.3.2]{virrion}).
Si $T$ est l'ensemble vide, nous omettrons de l'indiquer dans toutes ces opérations.
Sauf mention du contraire, on supposera (sans nuire à la généralité) les $k$-schémas réduits.

 Soient $\PP$ un $\V$-schéma formel lisse et $T$ un diviseur de $P$.
De façon similaire à \cite[4.2.2]{Beintro2},
on dispose du foncteur $\underset{\longrightarrow}{\lim}$ :
$\smash{\underset{^{\longrightarrow}}{LD}} ^{\mathrm{b}} _{\Q ,\mathrm{qc}}
(\smash{\widehat{\D}} _{\PP} ^{(\bullet)}(T))
\rightarrow
D  ( \smash{\D} ^\dag _{\PP} (\hdag T) _{\Q} )$. Celui-ci induit une équivalence
entre $D ^\mathrm{b} _\mathrm{coh} ( \smash{\D} ^\dag _{\PP} (\hdag T) _{\Q} )$
et une sous-catégorie pleine de $\smash{\underset{^{\longrightarrow}}{LD}} ^{\mathrm{b}} _{\Q ,\mathrm{qc}}
(\smash{\widehat{\D}} _{\PP} ^{(\bullet)}(T))$, notée
$\smash{\underset{^{\longrightarrow}}{LD}} ^{\mathrm{b}} _{\Q ,\mathrm{coh}}
(\smash{\widehat{\D}} _{\PP} ^{(\bullet)}(T))$ (voir \cite[4.2.4]{Beintro2} et \cite[1.1.3]{caro_courbe-nouveau}).
On dispose du même résultat pour les $F$-complexes.
Il nous arrivera d'omettre d'indiquer le foncteur $\underset{\longrightarrow}{\lim}$, surtout
lorsque l'on aura affaire à des $F$-complexes cohérents.

De plus, on désigne par $\sp $ : $\PP _K \rightarrow \PP$ le morphisme de spécialisation
du $\V$-schéma formel lisse $\PP$ ci-dessus (voir \cite[0.2.2.1]{Berig}).
Soient $X, Z$ deux sous-schémas fermés de $P$.
  On note $j _{Z,X}$ l'inclusion ouverte $X \setminus Z \subset X$ ou
  simplement $j _Z$ lorsque aucune confusion n'est à craindre.
  On dispose alors du foncteur $j ^\dag _Z$ (voir \cite[2.1.1]{Berig}) de la catégorie des $\O _{]X[ _{\PP}}$-modules
  dans celle des $j ^\dag _Z \O _{]X[ _{\PP}}$-modules
  (rappelons que $]X[ _{\PP}:= \sp ^{-1} (X)$ est le tube de $X$ dans $\PP _K$).
 La catégorie des $j _T ^\dag \O _{]X [ _{\PP}}$-modules cohérents munis d'une connexion surconvergente et d'une structure
  compatible de Frobenius sera notée
  $F\text{-}\mathrm{Isoc} ^{\dag} (\PP , T , X /K)$.
En notant $Y := X \setminus T$, cette catégorie est canoniquement isomorphe à celle des $F$-isocristaux sur $Y$ surconvergents le long de $T$ et notée
  $F\text{-}\mathrm{Isoc} ^{\dag} (Y ,X/K)$ (voir \cite[2.3.2 et 2.3.7]{Berig}). Ainsi, lorsque $P $ est propre,
  $F\text{-}\mathrm{Isoc} ^{\dag} (Y /K)\cong
  F\text{-}\mathrm{Isoc} ^{\dag} (\PP, T , X /K) $.
Si $\PP$ est le complété $p$-adique d'un $\V$-schéma formel faible propre et lisse $P ^\dag$,
on pose
$F\text{-}\mathrm{Isoc} ^{\dag} (P ^\dag , T , X /K):=F\text{-}\mathrm{Isoc} ^{\dag} (\PP , T , X /K)$.
On note $cv _{Y}$, le foncteur canonique de restriction $ (F\text{-})\mathrm{Isoc} ^\dag (Y/K) \rightarrow (F\text{-})\mathrm{Isoc}  (Y/K) $.
Si $E \in (F\text{-})\mathrm{Isoc} ^\dag (Y/K) $,
on écrit aussi $\widehat{E}: =cv _{Y}(E) $.

\section{\label{prenota-Fisocdagetc}$\D$-modules arithmétiques associés aux $F$-isocristaux surconvergents}

\subsection{$F$-isocristaux surconvergents sur les variétés lisses}

Soit $Y$ un $k$-schéma séparé, lisse.
En général, il n'existe pas de $\V$-schéma formel propre et lisse $\PP$,
de diviseur $T$ de $P$ tel que $Y$ soit un sous-schéma fermé de $\PP \setminus T$.
Cependant, nous donnons dans cette section une procédure de recollement
qui nous permet,
en nous ramenant à cette situation géométrique,
de construire la catégorie
$F \text{-} \mathrm{Isoc} ^{\dag} ( Y/K) $ des $F$-isocristaux surconvergents sur $Y$
(voir \ref{ind-isocsurcvgen-vide}).
Cette procédure est surtout intéressante car elle sera exactement reprise pour définir
les {\og $F$-isocristaux surcohérents sur $Y$\fg}
(voir \ref{ind-isosurcoh}).

\begin{nota}
\label{isocsurcvgen-vide}
Soit $Y$ un $k$-schéma séparé, lisse, intègre.
On choisit un recouvrement fini ouvert $(Y _\alpha) _{\alpha \in \Lambda}$ de $Y$ tel que,
pour tout $\alpha\in \Lambda$, il existe
un $\V$-schéma formel propre et lisse $\PP _\alpha$,
  un sous-schéma fermé $X _\alpha$ de
  $P _\alpha$ et un diviseur $T_\alpha$ de $P_\alpha$ tels que $Y _\alpha = X _\alpha\setminus T_\alpha$.
 De tels recouvrements existent bien. En effet,
 si $Y _\alpha$ est un ouvert affine de $Y$, il existe une immersion fermée de la forme
 $Y _\alpha \hookrightarrow \A ^n _k$. Il suffit alors de prendre $\PP _\alpha$ égal au
  complété $p$-adique de l'espace projectif $\P ^n _\V$, $T _\alpha$ le diviseur (réduit) de $\P ^n _k$
  complémentaire de $\A ^n _k$ et $X _\alpha$ l'adhérence schématique de $Y _\alpha $ dans $\P ^n _k$.

  Pour tous $\alpha, \beta$ et $\gamma \in \Lambda$,
  on note $p _1 ^{{\alpha \beta}}$ : $\PP _\alpha \times _\S \PP _\beta \rightarrow \PP _\alpha$ et
  $p _2 ^{{\alpha \beta}}$ : $\PP _\alpha \times _\S \PP _\beta \rightarrow \PP _\beta$
  les projections canoniques, $X _{{\alpha \beta}}$
  l'adhérence schématique de
  $Y _\alpha \cap Y _\beta$ dans $P _\alpha \times P _\beta$ (via l'immersion
  $Y _\alpha \cap Y _\beta \hookrightarrow Y _\alpha \times Y _\beta
  \hookrightarrow P _\alpha \times P _\beta$),
  $T _{{\alpha \beta}}=
  \smash{p _1 ^{\alpha \beta }} ^{-1} (T _\alpha) \cup \smash{p _2 ^{\alpha \beta }} ^{-1} (T _\beta)$.
  On remarque que
  $F\text{-}\mathrm{Isoc} ^{\dag} (\PP _\alpha \times \PP _\beta , T _{\alpha \beta}, X _{\alpha \beta}/K)
  \cong
  F\text{-}\mathrm{Isoc} ^{\dag} (Y _\alpha \cap Y _\beta /K )$.
  En notant $X _{\alpha \beta \gamma}$ l'adhérence schématique de
  $Y _\alpha \cap Y _\beta \cap Y _\gamma$ dans $P _\alpha \times P _\beta \times P _\gamma$ et
  $T _{\alpha \beta \gamma}$ la réunion des images inverses
  de $T _\alpha$, $T _\beta$, $T_\gamma$ par les projections de
  $\PP _\alpha \times \PP _\beta \times \PP _\gamma$ sur $\PP _\alpha$, $\PP _\beta$, $\PP _\gamma$, on obtient
  $F\text{-}\mathrm{Isoc} ^{\dag} (Y _\alpha \cap Y _\beta \cap Y _\gamma /K)\cong
  F\text{-}\mathrm{Isoc} ^{\dag}
  (\PP _\alpha \times \PP _\beta \times \PP _\gamma ,
  T _{\alpha \beta \gamma}, X _{\alpha \beta \gamma}/K)$.
On dispose des mêmes remarques en remplaçant
{\og $F\text{-}\mathrm{Isoc} ^{\dag}$\fg}
par {\og $F\text{-}\mathrm{Isoc} ^{\dag \dag}$\fg} (voir \cite[6.4]{caro_devissge_surcoh}).

  Soient
  $j _1 ^{\alpha \beta} $ : $Y _{\alpha}  \cap Y _{\beta} \subset Y _\alpha$,
  $j _2 ^{\alpha \beta} $ : $Y _{\alpha}  \cap Y _{\beta} \subset Y _\beta$,
  $j _{12} ^{\alpha \beta \gamma} $ : $Y _{\alpha}  \cap Y _{\beta}\cap Y _{\gamma} \subset
  Y _{\alpha}  \cap Y _{\beta}$,
  $j _{23} ^{\alpha \beta \gamma} $ :
  $Y _{\alpha}  \cap Y _{\beta}\cap Y _{\gamma}\subset Y _\beta \cap Y _\gamma$
  et
  $j _{13} ^{\alpha \beta \gamma} $ : $Y _{\alpha } \cap Y _{\beta}  \cap Y _{\gamma}
  \subset
  Y _\alpha \cap Y _\gamma$
  les immersions ouvertes.

La projection canonique d'espaces analytiques rigides
$p _{1K} ^{\alpha \beta}$ : $\PP _{\alpha K} \times \PP _{\beta K} \rightarrow \PP _{\alpha K}$
induit le morphisme ${{]X _{\alpha \beta} [ }_{\PP _\alpha \times \PP _\beta}}
\rightarrow {{]X _{\alpha} [ }_{\PP _\alpha}}$ que l'on note encore $p _{1K} ^{\alpha \beta}$.
D'où le foncteur canonique $j _1 ^{\alpha \beta *}:=
j ^\dag _{T _{\alpha \beta}} p _{1K} ^{\alpha \beta*}$ :
$F\text{-}\mathrm{Isoc} ^{\dag} ( \PP _\alpha, T _\alpha, X _\alpha /K)
\rightarrow F\text{-}\mathrm{Isoc} ^{\dag} ( \PP _{\alpha } \times \PP _\beta, T _{\alpha \beta}, X _{\alpha \beta} /K)$.
On définit de façon similaire $j _2 ^{\alpha \beta *}$, $j _{12} ^{\alpha \beta \gamma *} $,
$j _{13} ^{\alpha \beta \gamma*} $, $j _{23} ^{\alpha \beta \gamma*} $.

D'après \cite[6.4.3.b)]{caro_devissge_surcoh},
on bénéficie aussi du foncteur $j _1 ^{\alpha \beta *}:=
\R \underline{\Gamma} ^\dag _{X _{{\alpha \beta}}}
\circ (\hdag T _{{\alpha \beta}})\circ p _1 ^{\alpha \beta !}$ :
$F\text{-}\mathrm{Isoc} ^{\dag \dag} ( \PP _\alpha, T _\alpha, X _\alpha /K)
\rightarrow F\text{-}\mathrm{Isoc} ^{\dag \dag} ( \PP _\alpha \times \PP _\beta, T _{\alpha \beta}, X _{\alpha \beta} /K)$,
etc.

\end{nota}

\begin{defi}
\label{isocsurcvgen-defi}
\begin{enumerate}
  \item
Avec les notations \ref{isocsurcvgen-vide},
    on construit de la façon suivante la catégorie $F\text{-}\mathrm{Isoc} ^{\dag}
  (Y,\, (Y _\alpha, \PP _\alpha, T_\alpha,  X_\alpha) _{\alpha \in \Lambda}/K)$ :
\begin{itemize}
  \item
  Un objet est constitué par la donnée, pour tout $\alpha \in \Lambda$,
  d'un élément $E _\alpha$ de
  $F\text{-}\mathrm{Isoc} ^{\dag} ( \PP _\alpha, T _\alpha, X _\alpha /K)$ et,
  pour tous $\alpha,\beta \in \Lambda$,
  d'un isomorphisme
  $\eta _{\alpha \beta}$ : $j _2 ^{\alpha \beta *} (E _\beta)
  \riso j _1 ^{\alpha \beta *} (E _\alpha)$
  dans
    $F\text{-}\mathrm{Isoc} ^{\dag} (\PP _\alpha \times \PP _\beta , T _{\alpha \beta}, X _{\alpha \beta}/K)$,
ces isomorphismes vérifiant dans $F\text{-}\mathrm{Isoc} ^{\dag} (Y _\alpha \cap Y _\beta \cap Y _\gamma /K)$
la condition de cocycle
$j _{13}
^{\alpha \beta \gamma *} (\eta _{  \alpha \gamma} )=
j _{12} ^{\alpha \beta \gamma *} (\eta _{  \alpha \beta} )
\circ j _{23} ^{\alpha \beta \gamma *} ( \eta _{ \beta \gamma })$.
La famille d'isomorphismes $(\eta _{\alpha \beta}) _{\alpha ,\beta \in \Lambda}$
est appelée {\it donnée de recollement} de $ (E _\alpha) _{\alpha \in \Lambda} $.

\item Une flèche
$ ((E _\alpha) _{\alpha \in \Lambda}, (\eta _{  \alpha \beta}) _{\alpha,\beta  \in \Lambda}) \rightarrow
  ((E '_\alpha) _{\alpha \in \Lambda}, (\eta '_{  \alpha \beta}) _{\alpha,\beta  \in \Lambda})$
  est une famille de morphismes
  $ E _\alpha \rightarrow   E '_\alpha$
  de $F\text{-}\mathrm{Isoc} ^{\dag} ( \PP _\alpha, T _\alpha, X _\alpha /K)$
  commutant aux données de recollement respectives.

\end{itemize}
\item
En remplaçant les expressions
{\og $F\text{-}\mathrm{Isoc} ^{\dag}$\fg}
par {\og $F\text{-}\mathrm{Isoc} ^{\dag \dag}$\fg},
 on définit de même $F\text{-}\mathrm{Isoc} ^{\dag \dag}
  (Y,\, (Y _\alpha, \PP _\alpha, T_\alpha,  X_\alpha) _{\alpha \in \Lambda}/K)$.
Par convention, on notera plutôt dans ce cas
{\og $\E$ \fg} pour {\og $E$ \fg},
{\og $\theta _{\alpha \beta} $\fg} pour {\og $\eta _{\alpha \beta} $\fg}.
\end{enumerate}

Lorsque les diviseurs $T _\alpha$ sont vides, on écrit
$F\text{-}\mathrm{Isoc} ^{\dag}
  (X,\, (X _\alpha, \PP _\alpha) _{\alpha \in \Lambda}/K)$
  et
  $F\text{-}\mathrm{Isoc} ^{\dag \dag}
  (X,\, (X _\alpha, \PP _\alpha) _{\alpha \in \Lambda}/K)$.
\end{defi}

\begin{rema}
\label{rema-isocsurcvgen-vide}
  On garde les notations \ref{isocsurcvgen-defi} mais on remplace
  l'hypothèse {\og il existe
un $\V$-schéma formel propre et lisse $\PP _\alpha$\fg}
par
{\og il existe
un $\V$-schéma formel séparé et lisse $\PP _\alpha$\fg}.
  On définit alors
  Les catégories
  $F\text{-}\mathrm{Isoc} ^{\dag}
  (Y,\, (Y _\alpha, \PP _\alpha, T_\alpha,  X_\alpha) _{\alpha \in \Lambda}/K)$
  et
  $F\text{-}\mathrm{Isoc} ^{\dag\dag}
  (Y,\, (Y _\alpha, \PP _\alpha, T_\alpha,  X_\alpha) _{\alpha \in \Lambda}/K)$
  de la même façon qu'en \ref{isocsurcvgen-defi}.
\end{rema}

\begin{prop}
\label{ind-isocsurcvgen-vide}
On conserve les notations et hypothèses de \ref{isocsurcvgen-defi}.
La catégorie
$F\text{-}\mathrm{Isoc} ^{\dag}
  (Y,\, (Y _\alpha, \PP _\alpha, T_\alpha,  X_\alpha) _{\alpha \in \Lambda})$
  ne dépend canoniquement pas du choix de la famille
$(Y _\alpha, \PP _\alpha, T_\alpha,  X_\alpha) _{\alpha \in \Lambda}$.
  On la note alors sans ambiguïté
$F\text{-}\mathrm{Isoc} ^{\dag } (Y /K)$.

Si $Y$ est un $k$-schéma séparé, lisse et si
$(Y _{l}) _{l=1,\dots ,N}$ sont ses composantes connexes,
on pose alors
$F\text{-}\mathrm{Isoc} ^{\dag } (Y /K):= \prod _l F\text{-}\mathrm{Isoc} ^{\dag } (Y _l/K)$.
 Ses objets sont les
{\it $F$-isocristaux surconvergents sur $Y$}.
\end{prop}

\begin{proof}
Cela résulte des trois étapes suivantes :

i) Supposons d'abord qu'il existe un $\V$-schéma formel propre et lisse $\PP$, un sous-schéma fermé $X$ de $P$
tels que,
pour tout $\alpha \in \Lambda$, $\PP = \PP _\alpha$, $X = X _\alpha$ et
$T :=\cap _{\alpha\in \Lambda} T _\alpha$ soit un diviseur de $P$.
Prouvons alors que le foncteur canonique
\begin{equation}
  \label{ind-isocsurcvgen-videi)}
F\text{-}\mathrm{Isoc} ^{\dag} (\PP, \, T , \, X/K)
\rightarrow
F\text{-}\mathrm{Isoc} ^{\dag} ( Y,\, (Y _\alpha, \PP , T_\alpha,  X) _{\alpha \in \Lambda}/K),
\end{equation}
qui
à $E \in F\text{-}\mathrm{Isoc} ^{\dag} (\PP, \, T , \, X/K)$ associe
$((j^\dag _{T _\alpha} E) _{\alpha \in \Lambda}, (\eta _{  \alpha \beta}) _{\alpha,\beta  \in \Lambda})$,
où
$\eta _{\alpha \beta}$ :
$j ^\dag _{T _{\alpha \beta}} p _{2K} ^{\alpha \beta*} j^\dag _{T _\beta}   E
\riso
j ^\dag _{T _{\alpha \beta}} p _{1K} ^{\alpha \beta*}  j^\dag _{T _\alpha} E$
est l'isomorphisme canonique,
est une équivalence de catégories
(en notant $\delta$ : $\PP \hookrightarrow \PP \times \PP$ l'immersion fermée diagonale,
$\eta _{\alpha \beta}$ est l'unique morphisme induisant
l'isomorphisme canonique
$\delta ^* _K (\eta _{\alpha \beta})$ :
$j^\dag _{T _\alpha \cup T _\beta} j^\dag _{T _\beta}   E
\riso  j^\dag _{T _\alpha \cup T _\beta} j^\dag _{T _\alpha} E $).

D'abord, sa pleine fidélité est immédiate. \'Etablissons à présent son essentielle
surjectivité.
Or, avec les notations de \ref{rema-isocsurcvgen-vide}
et en posant $\PP _\alpha := \PP \setminus T _\alpha$ et
$X _\alpha := X \setminus T _\alpha$,
on bénéficie du foncteur canonique :
\newline \noindent
$F\text{-}\mathrm{Isoc} ^{\dag} ( Y,\, (Y _\alpha, \PP , T_\alpha,  X) _{\alpha \in \Lambda}/K)
\rightarrow
F\text{-}\mathrm{Isoc} ^{\dag} ( Y,\, (Y _\alpha, \PP _\alpha , T_\alpha,  X _\alpha) _{\alpha \in \Lambda}/K)$
défini par
$((E _\alpha) _{\alpha \in \Lambda}, (\eta _{  \alpha \beta}) _{\alpha,\beta  \in \Lambda})
\mapsto
(( E _\alpha |]X [ _{\PP} \cap \PP _{\alpha K}) _{\alpha \in \Lambda},
(\eta _{  \alpha \beta}
|]X [ _{\PP \times \PP} \cap \PP _{\alpha K} \times \PP _{\beta K}) _{\alpha,\beta  \in \Lambda})$.
On vérifie sans peine la pleine fidélité de ce dernier.
De plus, d'après la procédure de recollement de \cite[2.3.2.iii)]{Berig}, le foncteur composé obtenu
$F\text{-}\mathrm{Isoc} ^{\dag} (\PP, \, T , \, X/K)
\rightarrow
F\text{-}\mathrm{Isoc} ^{\dag} ( Y,\, (Y _\alpha, \PP _\alpha , T_\alpha,  X _\alpha) _{\alpha \in \Lambda}/K)$
est une équivalence de catégories. Il en résulte que
l'essentielle surjectivité du foncteur canonique
$F\text{-}\mathrm{Isoc} ^{\dag} (\PP, \, T , \, X/K)
\rightarrow
F\text{-}\mathrm{Isoc} ^{\dag} ( Y,\, (Y _\alpha, \PP , T_\alpha,  X) _{\alpha \in \Lambda}/K)$.

ii) Soit
un recouvrement fini ouvert $(Y _\alpha) _{\alpha \in \Lambda}$
(resp. $(Y '_{\alpha '}) _{\alpha '\in \Lambda '}$)
de $Y$ (resp. $Y'$) tel qu'il existe un $\V$-schéma formel propre et lisse $\PP _\alpha$ (resp. $\PP '_{\alpha'}$),
  un sous-schéma fermé $X _\alpha$ de $P _\alpha$ (resp $X '_{\alpha'}$ de $P '_{\alpha'}$)
  et un diviseur $T_\alpha$ de $P_\alpha$ (resp. $T '_{\alpha'}$ de $P '_{\alpha'}$)
  tels que
  $Y _\alpha = X _\alpha\setminus T_\alpha$
  (resp. $Y '_{\alpha'}= X '_{\alpha'} \setminus T '_{\alpha'}$).
  Grâce à l'étape i), quitte à prendre un recouvrement plus fin que
  $(Y _\alpha) _{\alpha \in \Lambda}$ et $(Y '_{\alpha '}) _{\alpha '\in \Lambda '}$,
  on se ramène à traiter le cas où $\Lambda =\Lambda '$ et $Y _\alpha = Y '_{\alpha}$.

iii) Quitte à remplacer $\PP ' _\alpha$ par $\PP' _\alpha \times _\S \PP _\alpha$,
on peut supposer qu'il existe un morphisme propre et lisse
$f _\alpha$ : $\PP ' _\alpha \rightarrow \PP _\alpha$
dont la restriction à $Y _\alpha$ induise l'identité de $Y _\alpha$.
Or,
on dispose de l'équivalence canonique de catégories
$Id _{Y _\alpha} ^* $ : $F\text{-}\mathrm{Isoc} ^{\dag} (\PP _\alpha, \,  X _\alpha, \, T_\alpha/K )
\cong
F\text{-}\mathrm{Isoc} ^{\dag} (\PP '_\alpha ,\,  X' _\alpha ,\, T '_\alpha/K)$ (voir \cite[2.3.5.i)]{Berig}).
Si
$ (E _\alpha) _{\alpha \in \Lambda}$ est une famille
d'objets $E _\alpha$ de $F\text{-}\mathrm{Isoc} ^{\dag} (\PP _\alpha, T _\alpha, X _\alpha /K)$
munie d'une donnée de recollement,
on lui associe la famille $ (Id ^* _{Y _\alpha}  (\E _\alpha) ) _{\alpha \in \Lambda}$, qui est munie d'une structure
canonique induite de donnée de recollement.
Ce foncteur est une équivalence de catégories.
\hfill \hfill \qed \end{proof}

\subsection{Construction : cas de la compactification lisse}

  Avec les notations et hypothèses de \ref{rema-isocsurcvgen-vide},
  on suppose de plus qu'il existe
  un $\V$-schéma formel séparé et lisse $\PP$, un sous-schéma fermé {\it lisse} intègre $X$ de $P$
tels que,
pour tout $\alpha \in \Lambda$, $\PP = \PP _\alpha$, $X = X _\alpha$,
 $T: =\cap _{\alpha\in \Lambda} T _\alpha$ soit un diviseur de $P$.
L'équivalence de catégories \ref{ind-isocsurcvgen-videi)}
reste encore valable dans ce cas.
On dispose ainsi d'une équivalence canonique de catégories
$F\text{-}\mathrm{Isoc} ^{\dag} ( Y,X/K)
\cong
F\text{-}\mathrm{Isoc} ^{\dag} ( Y,\, (Y _\alpha, \PP , T_\alpha,  X) _{\alpha \in \Lambda}/K)$.

Soit  $ ((E _\alpha) _{\alpha \in \Lambda}, (\eta _{  \alpha \beta}) _{\alpha,\beta  \in \Lambda})  \in
F\text{-}\mathrm{Isoc} ^{\dag} ( Y,\, (Y _\alpha, \PP , T_\alpha,  X) _{\alpha \in \Lambda}/K)$.
On pose $\E _\alpha := \sp _{X \hookrightarrow \PP , T _\alpha, +} (E _\alpha)$ (voir \cite{caro-construction})
et on définit l'isomorphisme $\theta _{  \alpha \beta}$ via la commutativité du diagramme :
\small
$$\xymatrix @R=0,3cm {
{\R \underline{\Gamma} ^\dag _{X _{{\alpha \beta}}}
\circ (\hdag T _{{\alpha \beta}})\circ p _2 ^{\alpha \beta !}
\sp _{X \hookrightarrow \PP , T _\beta, +} (E _\beta) }
\ar@{.>}[r] _-\sim ^-{\theta _{  \alpha \beta}}
&
{\R \underline{\Gamma} ^\dag _{X _{{\alpha \beta}}}
\circ (\hdag T _{{\alpha \beta}})\circ p _1 ^{\alpha \beta !}
(\sp _{X \hookrightarrow \PP , T _\alpha, +} (E _\alpha) )}
\\
{\sp _{X \hookrightarrow \PP\times \PP , T _{\alpha \beta}, +}
(j ^\dag _{T _{\alpha \beta}} p _{2K} ^{\alpha \beta*}    E_\beta)}
\ar[r] _-\sim ^-{\eta _{  \alpha \beta}} \ar[u] _-\sim
&
{\sp _{X \hookrightarrow \PP\times \PP , T _{\alpha \beta}, +}
(j ^\dag _{T _{\alpha \beta}} p _{1K} ^{\alpha \beta*} E _\alpha),}
\ar[u] _-\sim
}$$
\normalsize
où les flèches verticales viennent de \cite[4.1.2, 4.1.8]{caro-construction}.
D'où le foncteur
$\sp _+$ :
$F\text{-}\mathrm{Isoc} ^{\dag} ( Y,\, (Y _\alpha, \PP , T_\alpha,  X) _{\alpha \in \Lambda}/K)
\rightarrow
F\text{-}\mathrm{Isoc} ^{\dag \dag} ( Y,\, (Y _\alpha, \PP , T_\alpha,  X) _{\alpha \in \Lambda}/K)$
défini par
$((E _\alpha) _{\alpha \in \Lambda}, (\eta _{  \alpha \beta}) _{\alpha,\beta  \in \Lambda})
\mapsto
((\E _\alpha) _{\alpha \in \Lambda}, (\theta _{  \alpha \beta}) _{\alpha,\beta  \in \Lambda}) $.
Puisque notamment
$\sp _{X \hookrightarrow \PP , T _\alpha, +}$ :
$F\text{-}\mathrm{Isoc} ^{\dag} (\PP, \, T _\alpha , \, X /K)
\rightarrow
F\text{-}\mathrm{Isoc} ^{\dag\dag} (\PP, \, T _\alpha , \, X /K)$
est une équivalence de catégories (ainsi que
$\sp _{X \hookrightarrow \PP\times \PP , T _{\alpha \beta}, +}$, etc.),
on vérifie qu'il en est de même de
$\sp _+$.
Il en résulte par composition l'équivalence de catégories :
\begin{equation}
  \label{congcascompliss}
  F\text{-}\mathrm{Isoc} ^{\dag} ( Y,X/K)
\cong
F\text{-}\mathrm{Isoc} ^{\dag \dag} ( Y,\, (Y _\alpha, \PP , T_\alpha,  X) _{\alpha \in \Lambda}/K).
\end{equation}

\subsection{Surcohérence des $F$-isocristaux surconvergents : cas affine et lisse}
Pour mémoire, on dispose de la proposition suivante (voir \cite[6.3.1]{caro}).
\begin{prop}
\label{affdiv}
Soit $\widetilde{k}$ un corps quelconque,
$\smash{\widetilde{P}}$ un $\smash{\widetilde{k}}$-schéma propre et lisse, $\smash{\widetilde{T}}$ un sous-schéma fermé de $\smash{\widetilde{P}}$
de complémentaire $\smash{\widetilde{U}}$.
Si $\smash{\widetilde{U}}$ est affine et dense dans $\smash{\widetilde{P}}$ alors $\smash{\widetilde{U}}$ est le support d'un diviseur.
\end{prop}

Le lemme suivant fournit une extension de \cite[6.5.1]{caro_devissge_surcoh}.
\begin{lemm}
  \label{651}
  Soient $\PP$ un $\V$-schéma formel propre et lisse, $T$ un diviseur de $P$, $\U :=\PP \setminus T$,
  $Y$ un sous-schéma fermé intègre et lisse de $U$, $X$ l'adhérence de $Y$ dans $P$.
  Soit le diagramme commutatif
\begin{equation}
    \label{diag631dejong}
  \xymatrix @R=0,3cm {
  { Y '} \ar[r] \ar[d] ^b & {X'} \ar[r] ^{u'} \ar[d] ^a & {\PP'} \ar[d] ^f \\
  {Y} \ar[r] & {X } \ar[r] ^{u} & {\PP,}}
\end{equation}
  où $f$ est un morphisme propre et lisse de $\V$-schémas formels lisses,
  le carré de gauche est cartésien, $X'$ est un $k$-schéma lisse, $u'$ est une immersion fermée,
  $a$ est un morphisme
  projectif, surjectif, génériquement fini et étale tel que
  $a ^{-1} (T \cap X)$ soit un diviseur à croisement normaux de $X'$ (voir \cite[6.3.1]{caro_devissge_surcoh}).
Notons $T' :=f ^{-1} (T)$, $\U ':=\PP' \setminus T'$,
$g$ : $\U '\rightarrow \U $ le morphisme induit.

  Soient $\E _1, \E _2\in F\text{-}\mathrm{Isoc} ^{\dag \dag}( \PP,T , X/K)$,
$\phi$ :
$\R \underline{\Gamma} ^\dag _{X'} f ^! (\E _1) \rightarrow \R \underline{\Gamma} ^\dag _{X'} f ^! (\E _2)$
un morphisme de
$\mathrm{Isoc} ^{\dag \dag}( \PP',T',X' /K)$
et $\psi$ : $\E _1 |_{\U} \rightarrow \E _2 |_{\U}$
un morphisme de
$\mathrm{Isoc}^\dag  (\U,  Y /K)$.
Si $\phi$ et $\psi$ induise dans $\mathrm{Isoc}^\dag  (\U',  Y' /K)$
le même morphisme
$\R \underline{\Gamma} ^\dag _{Y'}g ^! (\E _1 |_{\U} )
\rightarrow \R \underline{\Gamma} ^\dag _{Y'} g ^! (\E _2 |_{\U} )$, alors
il existe un et un seul morphisme $\theta $ : $\E _1 \rightarrow \E _2$
de $\mathrm{Isoc} ^{\dag \dag}( \PP,T,X /K)$
induisant $\phi$ et $\psi$.
\end{lemm}

\begin{proof}
Pour $i =1,2$, en reprenant la preuve de \cite[6.3.1]{caro_devissge_surcoh}, nous construisons la suite de morphismes
$\E _i \rightarrow f _{T,+} \R \underline{\Gamma} ^\dag _{X'} f _T ^! (\E _i) \rightarrow \E _i$
dont le composé est un isomorphisme
(nous avons pour cela besoin de la structure de Frobenius de $\E _i$
car nous utilisons le théorème de pleine fidélité de Kedlaya \cite{kedlaya_full_faithfull}).
De manière analogue à la preuve de \cite[6.5.1]{caro_devissge_surcoh}, nous construisons alors
le morphisme $\theta $ : $\E _1 \rightarrow \E _2$ de $\mathrm{Isoc} ^{\dag \dag}( \PP,T,X /K)$ (via le diagramme
\cite[6.5.1.1]{caro_devissge_surcoh}).
En calquant la suite de la preuve de \cite[6.5.1]{caro_devissge_surcoh}, il suffit alors d'établir la fidélité du
foncteur
$\R \underline{\Gamma} ^\dag _{X'}  f _{T} ^! f _{T +}$ :
$\mathrm{Isoc} ^{\dag \dag}( \PP',T',X' /K)
\rightarrow
D ^\mathrm{b} _\mathrm{coh} ( \smash{\D} ^\dag _{\PP'} (\hdag T') _\Q)$.
Cela se vérifie comme dans \cite[6.5.1]{caro_devissge_surcoh}
en utilisant \cite[4.1.1]{tsumono} et \cite[4.3.12]{Be1} à la place de
\cite[4.1.2]{tsumono} et \cite{kedlaya_full_faithfull}.
\hfill \hfill \qed \end{proof}

\begin{lemm}
\label{TT'cas-lisse}
  Soient $\PP$ un $\V$-schéma formel séparé et lisse, $X$ un sous-schéma fermé lisse de $P$,
  $T\subset T'$ deux diviseurs de $P$ tels que $T \cap X$ et $T' \cap X$ soient des diviseurs de $X$,
  $E$ un ($F$-)isocristal sur $X \setminus T$ surconvergent le long de $T$.
  En notant
  $\E' : = \sp _{X \hookrightarrow \PP, T',+} ( j^\dag _{T'}E)$,
  on obtient
  $\E ' , \DD _{T} (\E ') \in (F\text{-})D ^\mathrm{b} _\mathrm{surcoh} ( \D ^\dag _{\PP} (\hdag T) _\Q)$.
\end{lemm}

\begin{proof}
  Cela découle de \cite[6.1.4]{caro_devissge_surcoh} et de
  $\E ' \riso (\hdag T') \sp _{X \hookrightarrow \PP, T,+} (E)$ (voir \cite[4.1.2]{caro-construction}).
\hfill \hfill \qed \end{proof}

\begin{lemm}
\label{plfdttssfrob}
  Avec les notations et hypothèses de \ref{TT'cas-lisse},
  le foncteur $(\hdag T')$ :
  $(F\text{-})\mathrm{Isoc} ^{\dag \dag} ( \PP ,\, T, \,X/K )
  \rightarrow
  (F\text{-})\mathrm{Isoc} ^{\dag \dag} ( \PP ,\, T', \,X/K )$ est pleinement fidèle.
\end{lemm}
\begin{proof}
  Cela résulte du diagramme commutatif à isomorphisme canonique près (voir \cite[4.1.2]{caro-construction}) :
  $$ \xymatrix @R=0,3cm  {
  {(F\text{-})\mathrm{Isoc} ^{\dag \dag} ( \PP ,\, T, \,X/K )}
  \ar[r] _-{(\hdag T')}
  &
  {(F\text{-})\mathrm{Isoc} ^{\dag \dag} ( \PP ,\, T', \,X/K )}
  \\
  {(F\text{-})\mathrm{Isoc} ^{\dag } ( Y/K)}
  \ar[r] _-{j ^\dag _{T'}}
  \ar[u] ^-{\sp _{X \hookrightarrow \PP, T,+}} _-\cong
  &
  {(F\text{-})\mathrm{Isoc} ^{\dag } ( Y '/K),}
  \ar[u] _-{\sp _{X \hookrightarrow \PP, T',+}} ^-\cong
  }$$
où les flèches verticales sont des équivalences de catégories
(voir \cite[6.2.2]{caro_devissge_surcoh} avec des structures de Frobenius,
mais cela reste de même valable sans structures de Frobenius, via
\cite[2.5.10]{caro-construction}, \cite[2.2.12]{caro_courbe-nouveau} et \cite[6.1.4]{caro_devissge_surcoh})
et celle du bas est pleinement fidèle
(voir \cite[4.1.1]{tsumono}).
\hfill \hfill \qed \end{proof}

\begin{prop}
\label{dag=dagdag-affine}
Soient $P ^\dag$ un $\V$-schéma formel faible propre et lisse,
$T$ un diviseur de $P$, $U ^\dag  := P ^\dag \setminus T$,
$Y ^\dag \hookrightarrow U ^\dag $
une immersion fermée de $\V$-schémas formels faibles lisses avec $Y ^\dag$ affine.
Soit $E \in F\text{-}\mathrm{Isoc} ^\dag (Y/K)$.

Le $\D ^\dag _{\PP}(\hdag T) _\Q$-module
$\sp _{Y ^\dag \hookrightarrow U ^\dag, T,+} (E)$
est alors muni d'une structure canonique de $F$-$\D ^\dag _{\PP} (\hdag T) _\Q$-module.
En outre, $\sp _{Y ^\dag \hookrightarrow U ^\dag, T,+} (E)\in F\text{-}\mathrm{Isoc} ^{\dag \dag} (Y/K)$.
\end{prop}

\begin{proof}
Il ne coûte rien de supposer $Y$ intègre.
Notons $X$ l'adhérence schématique de $Y$ dans $P$, $\E := \sp _{Y ^\dag \hookrightarrow U ^\dag, T,+} (E)$.
Soit $Z $ un sous-schéma fermé (réduit) de $X$.
Prouvons, par récurrence sur $(\dim Z, \mathrm{cmax} Z) $ (pour l'ordre lexicographique),
  où $\mathrm{cmax} Z$ désigne le nombre de composantes irréductibles de $Z$ de dimension $\dim Z$,
  l'assertion suivante :
  \small
\begin{equation}
  \label{dag=dagdag-affine-rec}
\text{{\og Si $\R \underline{\Gamma} ^\dag _Z \E
\in D ^\mathrm{b} _\mathrm{coh} ( \D ^\dag _{ \PP } (\hdag T) _\Q)$,
alors
$\R \underline{\Gamma} ^\dag _Z \E,
\DD _T  \R \underline{\Gamma} ^\dag _Z \E
\in F \text{-}D ^\mathrm{b} _\mathrm{surcoh} ( \D ^\dag _{ \PP } (\hdag T) _\Q)$.\fg}}
\end{equation}
\normalsize
  Traitons d'abord le cas où $\dim Z =0$. Il existe alors un $\V$-schéma formel faible lisse $Z ^\dag$ qui relève $Z$.
  Le cas où $Z \subset T$ est immédiat. De même, quitte à utilise l'hypothèse de récurrence, on se ramène au cas où $T \cap Z=\emptyset$.
  En notant $u$ : $ Z ^\dag \hookrightarrow Y ^\dag $ un relèvement de l'immersion fermée canonique $Z \hookrightarrow Y$,
  il résulte de \cite[5.2.6]{caro_devissge_surcoh} l'isomorphisme :
$\sp _{Z ^\dag  \hookrightarrow U ^\dag, T,+}
(u ^* E)[d _{Z/Y}]
\riso
\R \underline{\Gamma} ^\dag _Z \E $.
Comme $Z$ est lisse, $Z \subset Z $ se désingularise idéalement et donc, d'après \cite[7.2.4]{caro_devissge_surcoh},
$\sp _{Z ^\dag  \hookrightarrow U ^\dag, T,+} (u ^* E)$ est
muni d'une structure canonique de $F$-$\D ^\dag _{\PP} (\hdag T) _\Q$-module
et
$\R \underline{\Gamma} ^\dag _Z \E ,\, \DD _T  \R \underline{\Gamma} ^\dag _Z \E
\in F\text{-}D ^\mathrm{b} _\mathrm{surcoh} ( \D ^\dag _{ \PP } (\hdag T) _\Q)$.
  Supposons à présent $\dim Z >0$.
  Soient $Z _1 , \dots , Z _r$ les composantes irréductibles de $Z$ avec $\dim Z _1 = \dim Z$.
 Si $T \supset Z _1$ alors $\R \underline{\Gamma} ^\dag _{Z_1} \E=0$ et donc
$\R \underline{\Gamma} ^\dag _{Z_2 \cup \dots \cup Z _r} \E \riso
\R \underline{\Gamma} ^\dag _{Z} \E$
  (voir \cite[2.2.16]{caro_surcoherent}).
  Quitte à utiliser l'hypothèse de récurrence, on se ramène donc au cas où $T \not \supset Z_1$.
  Soit
  $x \in Z _1 \setminus (Z _2 \cup \cdots \cup Z _r \cup T)$.
  Il existe un ouvert affine $\widetilde{U}$ de $P$ tel que $x \in \widetilde{U}$ et
  $\widetilde{U} \cap (Z _2 \cup \cdots \cup Z _r\cup T) = \emptyset$. Ainsi,
 $\widetilde{U} \cap Z =\widetilde{U} \cap Z _1$ et $\widetilde{U} \cap Z $ est dense dans $Z _1$.
  Par \ref{affdiv}, $\widetilde{T}:= P \setminus \widetilde{U}$ est le support d'un diviseur.
  Quitte à rétrécir $\widetilde{U}$, grâce à \cite[7.4.1]{caro_devissge_surcoh},
  on peut en outre supposer
  $\widetilde{U} \cap Z$ lisse et muni d'un modèle idéal.
En notant $\smash{\widetilde{Y}}^\dag $ (resp. $\smash{\widetilde{U}}^\dag $)
l'ouvert de $Y ^\dag $ (resp. $U ^\dag $) complémentaire de $\widetilde{T}$
et $\widetilde{E}:= j _{\widetilde{T} } ^\dag E$
le $F$-isocristal surconvergent sur $\widetilde{Y}$ déduit de $E$,
on obtient $(\hdag \widetilde{T}) \E  \riso
\sp _{\widetilde{Y} ^\dag \hookrightarrow \widetilde{U} ^\dag, \widetilde{T},+} (\widetilde{E})$
(voir \cite[5.2.3]{caro_devissge_surcoh}).
Comme $\widetilde{Y}\cap Z=\widetilde{U}\cap Z$ est affine et lisse, il existe d'après Elkik (voir \cite{elkik})
un $\V$-schéma formel faible affine et lisse $\smash{\widetilde{Y}} ^\dag _Z$ le relevant.
Par \cite[5.2.6]{caro_devissge_surcoh}, comme
$\R \underline{\Gamma} ^\dag _Z (\hdag \widetilde{T}) \E
\riso (\hdag \widetilde{T})  \R \underline{\Gamma} ^\dag _Z \E
\in D ^\mathrm{b} _\mathrm{coh} ( \D ^\dag _{ \PP } (\hdag \widetilde{T}) _\Q)$,
en désignant par $a$ :
$\smash{\widetilde{Y}} ^\dag _Z \hookrightarrow \smash{\widetilde{Y}} ^\dag $
une immersion fermée relevant
$\widetilde{Y}\cap Z \hookrightarrow \widetilde{Y}$, il en résulte l'isomorphisme :
$\sp _{\smash{\widetilde{Y}} ^\dag _Z \hookrightarrow \widetilde{U} ^\dag, \widetilde{T},+}
(a ^* \widetilde{E})
\riso
\R \underline{\Gamma} ^\dag _Z (\hdag \widetilde{T}) \E [-d _{\widetilde{Y} _Z/\widetilde{Y}}]=:\FF.$
Avec \cite[7.2.4]{caro_devissge_surcoh}, il en découle que
$\FF
\in
F \text{-}\mathrm{Isoc}^{\dag \dag} (\PP,\widetilde{T},Z _1/K)$.
En utilisant le théorème de désingularisation de de Jong (voir \cite{dejong}), de manière analogue à
\cite[6.3.1]{caro_devissge_surcoh}, on obtient le diagramme commutatif :
$$\xymatrix @R=0,3cm {
{\widetilde{U}}
\ar[r]
&
{U}
\ar[r]
&
{\PP}
\ar@{=}[r]
&
{\PP}
\\
{\widetilde{Y}}
\ar[r]\ar[u]
&
{Y}
\ar[r]\ar[u]
&
{X}
\ar[r] \ar[u]
&
{\PP}
\ar@{=}[u]
\\
{\widetilde{Y}\cap Z _1}
\ar[r]\ar[u] ^-a
&
{Y\cap Z _1}
\ar[r]\ar[u]
&
{Z _1}
\ar[r]\ar[u]
&
{\PP}
\ar@{=}[u]
\\
{\widetilde{V} '}
\ar[r] \ar[u] ^-b
&
{V'}
\ar[r] \ar[u]
&
{Z'}
\ar[r] \ar[u] ^-h
&
{\PP ',}
\ar[u] ^-f
}$$
dans lequel les carrés des deux colonnes de gauche sont cartésiens,
$Z'$ est lisse,
$f$ est un morphisme propre et lisse de $\V$-schémas formels propres et lisses,
$h$ est un morphisme projectif, surjectif, génériquement fini et étale,
tel que $h ^{-1} ( Z _1 \cap \widetilde{T}) $ soit un diviseur de $Z'$,
les flèches horizontales des carrés de droite sont des immersions fermées.
On pose $T' := f ^{-1} (T)$ et $\widetilde{T}' := f ^{-1} (\widetilde{T})$.
Comme $Z'$ est lisse et $Z ' \cap \widetilde{T}'$ est un diviseur de $Z'$,
alors $Z'\cap T$ est aussi un diviseur de $Z'$.
D'après \cite[6.3.1]{caro_devissge_surcoh},
$\FF$
est un facteur direct de
  $f _{+} \R \underline{\Gamma} ^\dag _{Z'} f ^!
  (\FF)$
  et
il existe un
(unique à isomorphisme près)
$F$-isocristal $\smash{\widetilde{E}} '$ surconvergent sur $\smash{\widetilde{V}} '$
tel que
  $\R \underline{\Gamma} ^\dag _{Z'} f ^!
  (\FF)
  \riso
  \sp _{Z' \hookrightarrow \PP',\smash{\widetilde{T}}'} (\smash{\widetilde{E}} ')$.
Grâce au théorème de pleine fidélité de Kedlaya (voir \cite{kedlaya_full_faithfull}),
en notant $c$ : $\smash{\widetilde{V}} ' \rightarrow Y$ le morphisme canonique induit par $h$,
on vérifie $\smash{\widetilde{E}} ' \riso b ^* (a ^* \smash{\widetilde{E}}) \riso c ^* (E) $.
En particulier, $\smash{\widetilde{E}} '$ provient d'un $F$-isocristal surconvergent sur $V'$.
Par \ref{TT'cas-lisse},
$\R \underline{\Gamma} ^\dag _{Z'} f ^!
  (\FF)$
  et
  $\DD _{T'} \R \underline{\Gamma} ^\dag _{Z'} f ^!
  (\FF)\in F \text{-}D ^\mathrm{b} _\mathrm{surcoh} ( \D ^\dag _{\PP'} (\hdag T') _\Q)$.
  Par \cite[3.1.9]{caro_surcoherent} et \cite[1.2.7]{caro_courbe-nouveau},
  $f _+ \R \underline{\Gamma} ^\dag _{Z'} f ^!
  (\FF)$
  et
  $\DD _{T} f _+ \R \underline{\Gamma} ^\dag _{Z'} f ^!
  (\FF)\in
F \text{-}  D ^\mathrm{b} _\mathrm{surcoh} ( \D ^\dag _{\PP} (\hdag T) _\Q)$.
D'où
$\FF,
\DD _T \FF
\in
F \text{-}D ^\mathrm{b} _\mathrm{surcoh} ( \D ^\dag _{\PP} (\hdag T) _\Q)$.
Le triangle de localisation en $\widetilde{T}$ de
$\R \underline{\Gamma} ^\dag _Z \E $
s'écrit :
\begin{equation}
  \label{lemm-dag=dagdag-affine-trloc}
\R \underline{\Gamma} ^\dag _{Z \cap \widetilde{T}} \E
\rightarrow  \R \underline{\Gamma} ^\dag _Z \E \rightarrow
\R \underline{\Gamma} ^\dag _Z (\hdag \widetilde{T}) \E  \rightarrow +1.
\end{equation}
Il en dérive que
$\R \underline{\Gamma} ^\dag _{Z \cap \widetilde{T}} \E
\in D ^\mathrm{b} _\mathrm{coh} ( \D ^\dag _{\PP} (\hdag T) _\Q)$.
Comme $Z _1 \cap \widetilde{U}$ est dense dans $Z _1$, $\dim Z _1 \cap \widetilde{T} < \dim Z _1 = \dim Z$.
Ainsi, soit $\dim Z  \cap \widetilde{T} < \dim Z $ soit
$\dim Z  \cap \widetilde{T} = \dim Z $ et $\mathrm{cmac } Z  \cap \widetilde{T} < \mathrm{cmac } Z$.
Par hypothèse de récurrence, il en découle que
$\R \underline{\Gamma} ^\dag _{Z \cap \widetilde{T}} \E, \DD _T \R \underline{\Gamma} ^\dag _{Z \cap \widetilde{T}} \E
\in F \text{-}D ^\mathrm{b} _\mathrm{surcoh} ( \D ^\dag _{\PP} (\hdag T) _\Q)$.
On conclut la récurrence via
\ref{lemm-dag=dagdag-affine-trloc} et le triangle distingué induit en lui appliquant $\DD _T$.

En appliquant \ref{dag=dagdag-affine-rec} au cas où $Z =X$, on obtient
$\E \in F\text{-}\mathrm{Isoc} ^{\dag \dag} (Y/K)$.
Il reste à établir que $\E$ est muni d'une structure {\it canonique} de Frobenius.
Mise à part $X$ et $\E := \sp _{Y ^\dag \hookrightarrow U ^\dag, T,+} (E)$,
oublions à présent les notations introduites au cours de la preuve.

Grâce au théorème de désingularisation de de Jong (\cite{dejong}),
il existe un morphisme projectif et surjectif $a$ : $X' \rightarrow X $, qui se décompose
  en une immersion fermée $X ' \hookrightarrow \P ^r _{X }$ suivie de la projection
  canonique $\P ^r _{X } \rightarrow X $, tel que
    $X ' $ soit intègre et lisse, $a$ soit génériquement fini et étale
    et $a ^{-1} ( T \cap X)$ soit un diviseur de $X'$.
    On pose $P ^{\prime \dag} : =\P ^{r \dag} _{P ^\dag }$,
$f$ : $P ^{\prime \dag} \rightarrow P ^{\dag}$ la projection canonique,
$T' := f ^{-1}(T)$,
    $U ^{\prime \dag}:= P ^{\prime \dag} \setminus T '$,
$g$ : $U ^{\prime \dag} \rightarrow U ^{\dag}$ le morphisme induit par $f$.
La fin de la preuve s'établit via les trois étapes suivantes.

$(a)$ Si $T _1$ est un diviseur de $P$ contenant $T$, on pose alors
$Y ^\dag _1 := Y ^\dag \setminus T _1$, $Y ' _1  := a ^{-1} (Y _1)$, $T' _1 :=f ^{-1} (T _1)$
$U _1 ^{\dag}:= P ^\dag \setminus T _1$,
$U _1 ^{\prime \dag}:= P ^{\prime \dag} \setminus T' _1 $,
$g _1$ : $U _1 ^{\prime \dag} \rightarrow U _1 ^{\dag}$ le morphisme induit par $f$.
D'après \cite[7.4.1]{caro_devissge_surcoh} et sa preuve, on peut choisir $T _1$ tel que
$Y _1 := Y \setminus T _1$ soit affine et lisse,
le morphisme $b _1 $ : $Y ' _1 \rightarrow Y _1$ induit par $a$ soit fini et étale,
l'immersion fermée
$ Y _1 ' \hookrightarrow \P ^r _{Y _1}$ se relève en une immersion fermée de $\V$-schémas formels faibles lisses
de la forme
$Y _1 ^{\prime \dag} \hookrightarrow \P ^{r\dag} _{Y ^\dag _1}$ (en particulier $Y _1 \subset X$ se désingularise idéalement).
Notons $E _1 := j ^\dag _{T _1} (E)$.
On dispose des isomorphismes canoniques :
\begin{gather}\notag
\xymatrix @C=2cm {
  {\sp _{Y _1 ^{\prime \dag} \hookrightarrow U _1 ^{\prime \dag}, T' _1,+} ( b _1 ^* E _1)}
  \ar[r] _-\sim 
  &
{\R \underline{\Gamma} ^\dag _{X'} f ^!\sp _{Y _1 ^\dag \hookrightarrow U _1 ^\dag, T_1,+} (E _1)}
  \ar[r] _-\sim 
  &
  {}
  }
  \\
  \label{Y1Yiso}
\xymatrix @C=2cm {
    {}
  \ar[r] _-\sim 
  &
  {}
  }
(\hdag T'_1) \R \underline{\Gamma} ^\dag _{X'} f ^! \sp _{Y ^\dag \hookrightarrow U ^\dag, T,+} (E )=
(\hdag T'_1) \R \underline{\Gamma} ^\dag _{X'} f ^! \E .
\end{gather}
Comme $X'$ est lisse, $Y ' _1 \subset X'$ se désingularise idéalement.
Par \cite[7.2.4]{caro_devissge_surcoh}, le $\D ^\dag _{\PP'}(\hdag T '_1) _\Q$-module
$\sp _{Y _1 ^{\prime \dag} \hookrightarrow U _1 ^{\prime \dag}, T' _1,+} ( b _1 ^* E _1)$
est muni d'une structure canonique de Frobenius
(i.e. de $F \text{-}\D ^\dag _{\PP'}(\hdag T '_1) _\Q$-module) induite par celle de $E$.
Comme Frobenius commute canoniquement à toutes les opérations cohomologiques de la théorie des $\D$-modules arithmétiques
(voir \cite{Beintro2}, \cite{Be2}, \cite{caro_surcoherent} etc.),
on déduit alors de \ref{Y1Yiso}, l'isomorphisme canonique
$(\hdag T'_1) \R \underline{\Gamma} ^\dag _{X'} f ^! \E \riso (\hdag T'_1) \R \underline{\Gamma} ^\dag _{X'} f ^! F ^* \E $.
Or, comme $X'$ est lisse et que $T' _1\cap X'$, $T' \cap X'$ sont des diviseurs de $X'$,
le foncteur
$(\hdag T'_1)$ :
$\mathrm{Isoc} ^{\dag \dag} ( \PP ', T', X ' /K)\rightarrow \mathrm{Isoc} ^{\dag \dag} ( \PP ', T'_1, X ' /K)$
est, d'après \ref{plfdttssfrob}, pleinement fidèle.
Puisque $\R \underline{\Gamma} ^\dag _{X'} f ^!\E
\in \mathrm{Isoc} ^{\dag \dag} ( \PP ', T' , X ' /K)$, de même en remplaçant $\E$ par $F ^* \E$,
il en découle l'isomorphisme canonique
$\phi$ :
$\R \underline{\Gamma} ^\dag _{X'} f ^! (\E) \riso \R \underline{\Gamma} ^\dag _{X'} f ^! (F ^* \E )$.

$(b)$
On construit l'isomorphisme canonique :
$\psi$ : $\E  |_{\U} \rightarrow F ^* \E  |_{\U}$ via le diagramme commutatif suivant
\begin{equation}
  \label{dag=dagdag-affine-forbu}
  \xymatrix @R=0,4cm @C=2cm{
  {\sp _{Y \hookrightarrow \U, +} (\widehat{E})}
  \ar[d] _-\sim
  \ar@{=}[r]
  &
  {\sp _{Y \hookrightarrow \U, +} (\widehat{E})}
  \ar[r] _-\sim 
  \ar[d] _-\sim
  &
  {\E |\U}
  \ar@{.>}[d]_-\sim ^-{\psi}
  \\
{\sp _{Y \hookrightarrow \U, +} (F ^*\widehat{E})}
  \ar[r] _-\sim 
  &
  {F ^* \sp _{Y \hookrightarrow \U, +} (\widehat{E})}
  \ar[r] _-\sim 
  &
  {F ^* \E |\U,}
  }
\end{equation}
où la flèche de gauche provient de l'isomorphisme tautologique $\widehat{E} \riso F ^* \widehat{E}$
munissant $\widehat{E}$ d'une structure de Frobenius.

$(c)$ On vérifie par construction la commutativité du diagramme :
\begin{equation}
\label{U'U'1fid}
  \xymatrix @R=0,4cm @C=2cm{
  {\R \underline{\Gamma} ^\dag _{X'} f ^! (\E)| \U ' _1}
  \ar[r] _-\sim
  \ar[d] _-\sim ^-{\phi}
  &
  {\R \underline{\Gamma} ^\dag _{Y'} g  ^! (\E |\U) | \U ' _1}
  \ar[d]_-\sim ^-{\psi}
  \\
  {\R \underline{\Gamma} ^\dag _{X'} f ^! (F ^* \E ) | \U ' _1}
  \ar[r] ^-\sim
  &
  {\R \underline{\Gamma} ^\dag _{Y'} g  ^! (F ^* \E |\U)  | \U ' _1 ,}
  }
\end{equation}
Par \ref{plfdttssfrob} et \cite[4.3.12]{Be1},
$|\U ' _1$ : $\mathrm{Isoc}^\dag  (\U',  Y' /K)\rightarrow \mathrm{Isoc}^\dag  (\U'_1,  Y'_1 /K)$ est fidèle.
Le diagramme \ref{U'U'1fid} reste ainsi commutatif si on enlève les termes
{\og $| \U ' _1 $\fg}. Les isomorphismes $\phi$ et $\psi$ induisent donc canoniquement le même isomorphisme
de la forme
$\R \underline{\Gamma} ^\dag _{Y'} g  ^! (\E |\U) \riso
\R \underline{\Gamma} ^\dag _{Y'} g  ^! (F ^*\E |\U)$.
Il résulte alors de \ref{651} l'existence d'un isomorphisme canonique
$\E \riso F ^* \E $ induisant $\phi $ et $\psi$.
\hfill \hfill \qed \end{proof}

\begin{prop}
\label{eqcataffine}
  Soit $Y$ un $k$-schéma affine et lisse. On dispose d'une équivalence canonique de catégories
  $$\sp _{Y,+}\ : \ F\text{-}\mathrm{Isoc} ^\dag (Y/K) \cong F\text{-}\mathrm{Isoc} ^{\dag \dag} (Y/K) .$$
\end{prop}
\begin{proof}
Grâce au théorème de relèvement d'Elkik (voir \cite{elkik}),
  il existe $P ^\dag$ un $\V$-schéma formel faible propre et lisse,
$T$ un diviseur de $P$, $U ^\dag  := P ^\dag \setminus T$,
$Y ^\dag \hookrightarrow U ^\dag $
une immersion fermée de $\V$-schémas formels faibles lisses.
Par \cite[6.5.3]{caro_surcoherent}, nous disposons
  du foncteur pleinement fidèle $\rho _Y$ : $F\text{-}\mathrm{Isoc} ^{\dag \dag} (Y/K)
  \rightarrow F\text{-}\mathrm{Isoc} ^{\dag} (Y/K)$.
D'après \ref{dag=dagdag-affine},
on bénéficie de la factorisation :
$\sp _{Y ^\dag \hookrightarrow U ^\dag, T,+}$ :
$F\text{-}\mathrm{Isoc} ^{\dag} (Y/K) \rightarrow
F\text{-}\mathrm{Isoc} ^{\dag \dag} (Y/K)$.
En utilisant notamment \ref{dag=dagdag-affine-forbu} et
la peine fidélité du foncteur canonique de restriction
$F\text{-}\mathrm{Isoc} ^{\dag} (Y/K)
\rightarrow F\text{-}\mathrm{Isoc} (Y/K)$
(voir \cite{kedlaya_full_faithfull}),
on vérifie par construction que ces deux foncteurs, $\rho _Y$ et $\sp _{Y ^\dag \hookrightarrow U ^\dag, T,+}$,
sont quasi-inverses et induisent donc des équivalences de catégories.
\hfill \hfill \qed \end{proof}

\subsection{Construction : cas du recouvrement affine et lisse}
\begin{vide}
\label{isocrecolaff}
 Soit
$b$ : $Y ^{\prime \dag} \rightarrow Y ^{\dag}$ un morphisme de $\V$-schémas formels faibles affines et lisses.
La catégorie $F\text{-}\mathrm{Isoc} ^\dag (Y/K)$
est canoniquement isomorphe à celle des
$F$-$\D _{Y ^\dag, \Q }$-modules à gauche globalement de présentation finie, $\O _{Y ^\dag ,\Q }$-cohérent
(voir \cite[2.5.7]{Berig}).
Par \cite[2.5.6]{Berig}, le foncteur $b ^* = \O _{Y ^{\prime \dag} ,\Q } \otimes _{b  ^{-1} \O _{Y ^\dag ,\Q }}-$
est canoniquement isomorphe au foncteur de \cite[2.3.6]{Berig}
$b _0  ^* $ : $F\text{-}\mathrm{Isoc} ^\dag (Y/K) \rightarrow F\text{-}\mathrm{Isoc} ^\dag (Y'/K)$.
De plus, le foncteur canonique de restriction
$cv _{Y}$ : $ F\text{-}\mathrm{Isoc} ^\dag (Y/K) \rightarrow F\text{-}\mathrm{Isoc}  (Y/K) $
est canoniquement isomorphe à
$\O _{\Y ,\Q } \otimes _{\O _{Y ^\dag ,\Q }}-$ (voir \cite[5.1.3]{caro_devissge_surcoh}).
Le foncteur $\smash{\widehat{b}} ^* = \O _{\Y ^{\prime } ,\Q } \otimes _{b  ^{-1} \O _{\Y ,\Q }}-$
s'identifie canoniquement au foncteur image inverse
$F\text{-}\mathrm{Isoc}  (Y/K) \rightarrow F\text{-}\mathrm{Isoc}  (Y'/K)$
de \cite[2.3.2.2]{Berig}.
Enfin, on dispose d'un isomorphisme canonique
$cv _{Y '} \circ b ^* \riso \widehat{b} ^* \circ cv _{Y }$.
\end{vide}

\begin{vide}
\label{isocrecolaff2}
Avec les notations de \ref{isocrecolaff},
soient deux morphismes
$b,\, b'$ : $Y ^{\prime \dag} \rightarrow Y ^{\dag}$ de $\V$-schémas formels faibles affines et lisses
tels que $b _0 = b ' _0$.
Pour tout $E \in F\text{-}\mathrm{Isoc} ^\dag (Y/K)$,
on bénéficie d'un isomorphisme canonique
$b ^* (E)  \riso b ^{\prime *} (E)$,
transitif pour la composition de morphismes de $\V$-schémas formels faibles affines et lisses.
En effet,
on définit d'abord l'isomorphisme canonique
$cv _{Y '} (b ^* (E) ) \riso cv _{Y '} (b ^{\prime *} (E))$
via la commutativité du diagramme ci-après :
\begin{equation}
  \notag
  \xymatrix @R=0,3cm {
  {cv _{Y '} (b ^* (E) ) }
  \ar@{.>}[r]
  \ar[d] _-\sim
  &
  {cv _{Y '} (b ^{\prime *} (E))}
  \ar[d] _-\sim
  \\
  {\widehat{b} ^* cv _{Y }(E)}
  \ar[r] _-\sim
  &
  {\widehat{b} ^{\prime *} cv _{Y }(E),}
  }
\end{equation}
dont l'isomorphisme du bas est \cite[2.2.17.i)]{Berig}.
On conclut ensuite par pleine fidélité de $cv _{Y '} $ (voir \cite{kedlaya_full_faithfull}).
\end{vide}

\begin{vide}
\label{nota-Fisocdagetc}

Soit $Y$ un $k$-schéma séparé, lisse, intègre.
Soit $(Y _\alpha) _{\alpha \in \Lambda}$ un recouvrement fini de $Y$ par des ouverts affines et,
pour tout $\alpha\in \Lambda$, soient $P ^\dag _\alpha$ un $\V$-schéma formel faible propre et lisse,
  $X _\alpha$ un sous-schéma fermé de $P _\alpha$
  et $T_\alpha$ un diviseur de $P_\alpha$ tel que $Y _\alpha = X _\alpha\setminus T_\alpha$
  (on peut toujours en construire grâce à \cite{elkik}).

  Pour tous $\alpha, \beta, \gamma \in \Lambda$,
  on note $p _1 ^{{\alpha \beta}}$ : $P ^\dag _\alpha \times _\S P ^\dag _\beta \rightarrow P ^\dag _\alpha$ et
  $p _2 ^{{\alpha \beta}}$ : $P ^\dag _\alpha \times _\S P ^\dag _\beta \rightarrow P ^\dag _\beta$
  les projections canoniques, $X _{{\alpha \beta}}$
  l'adhérence schématique de
  $Y _\alpha \cap Y _\beta$ dans $P _\alpha \times P _\beta$,
  $T _{{\alpha \beta}}=
  \smash{p _1 ^{\alpha \beta }} ^{-1} (T _\alpha) \cup \smash{p _2 ^{\alpha \beta }} ^{-1} (T _\beta)$.
  De même, on note $X _{\alpha \beta \gamma}$ l'adhérence schématique de
  $Y _\alpha \cap Y _\beta \cap Y _\gamma$ dans $P _\alpha \times P _\beta \times P _\gamma$ et
  $T _{\alpha \beta \gamma}$ la réunion des images inverses
  de $T _\alpha$, $T _\beta$ et $T_\gamma$ par les projections de
  $P ^\dag _\alpha \times P ^\dag _\beta \times P ^\dag _\gamma$ sur $P ^\dag _\alpha$, $P ^\dag _\beta$
  et $P ^\dag _\gamma$.

Pour tous $\alpha ,\beta, \gamma \in \Lambda$, soient $Y ^\dag _{\alpha \beta}  $ et
$Y ^\dag _{\alpha \beta \gamma }$ des $\V$-schémas formels
faibles affines et lisses relevant respectivement $Y _{\alpha}  \cap Y _{\beta} $ et
$Y _{\alpha}  \cap Y _{\beta}\cap Y _{\gamma}$, soient
  $j _1 ^{\alpha \beta} $ : $Y ^\dag _{\alpha \beta}  \rightarrow Y ^\dag _\alpha$,
  $j _2 ^{\alpha \beta} $ : $Y ^\dag _{\alpha \beta}  \rightarrow Y ^\dag _\beta$,
  $j _{12} ^{\alpha \beta \gamma} $ :
  $Y ^\dag _{\alpha \beta \gamma}   \rightarrow Y ^\dag _{\alpha \beta} $,
  $j _{23} ^{\alpha \beta \gamma} $ :
  $Y ^\dag _{\alpha \beta \gamma}
  \rightarrow Y ^\dag _{\beta \gamma}   $
  et
  $j _{13} ^{\alpha \beta \gamma} $ :
  $Y ^\dag _{\alpha \beta \gamma}
\rightarrow
  Y ^\dag _{\alpha \gamma}   $
  des relèvements des immersions canoniques ouvertes.
  \end{vide}

\begin{vide}
Avec les notations et hypothèses de \ref{nota-Fisocdagetc},
on définit la catégorie
$F\text{-}\mathrm{Isoc} ^{\dag} ( Y,\, (Y _\alpha, P ^\dag _\alpha, T_\alpha,  X_\alpha) _{\alpha \in \Lambda}/K)$
dont un objet est constitué par la donnée, pour tout $\alpha \in \Lambda$,
  d'un élément $E _\alpha$ de
  $F\text{-}\mathrm{Isoc} ^{\dag } ( Y _\alpha /K)$ et,
  pour tous $\alpha,\beta \in \Lambda$,
  d'un isomorphisme dans
    $F\text{-}\mathrm{Isoc} ^{\dag} ( Y _{\alpha } \cap Y _{\beta} /K)$
de la forme
$j _2 ^{\alpha \beta *} (E _\beta)
  \riso j _1 ^{\alpha \beta *} (E _\alpha)$
  (voir \ref{isocrecolaff} à propos des foncteurs $j _1 ^{\alpha \beta *}$ et $j _2 ^{\alpha \beta *}$),
  ces isomorphismes vérifiant dans $F\text{-}\mathrm{Isoc} ^{\dag } (Y _\alpha \cap Y _\beta \cap Y _\gamma /K)$
la condition de cocycle
$j _{13}
^{\alpha \beta \gamma *} (\theta _{  \alpha \gamma} )=
j _{12} ^{\alpha \beta \gamma *} (\theta _{  \alpha \beta} )
\circ j _{23} ^{\alpha \beta \gamma *} ( \theta _{ \beta \gamma })$
(modulo les isomorphismes canoniques de \ref{isocrecolaff2}).

D'après \cite[2.5.6]{Berig},
le foncteur $j _1 ^{\alpha \beta *}$
(voir \ref{isocrecolaff}) est canoniquement isomorphe à celui défini dans
\ref{isocsurcvgen-vide}.
Il en découle l'équivalence canonique de catégories :
\begin{equation}
  \label{dag=dagafflisse}
 F\text{-}\mathrm{Isoc} ^{\dag}  (Y,\, (Y _\alpha, \PP _\alpha, T_\alpha,  X_\alpha) _{\alpha \in \Lambda}/K)
 \cong
F\text{-}\mathrm{Isoc} ^{\dag}  (Y,\, (Y _\alpha, P ^\dag _\alpha, T_\alpha,  X_\alpha) _{\alpha \in \Lambda}/K).
\end{equation}

\end{vide}

\begin{prop}
\label{recollcompsp+}
  Avec les notations et hypothèses \ref{nota-Fisocdagetc},
  on dispose d'une équivalence canonique de catégories :
  \begin{equation}
  \label{congcasaffliss}
  F\text{-}\mathrm{Isoc} ^{\dag}  (Y/K)
\cong
F\text{-}\mathrm{Isoc} ^{\dag \dag}  (Y,\, (Y _\alpha, \PP _\alpha, T_\alpha,  X_\alpha) _{\alpha \in \Lambda}/K).
\end{equation}
\end{prop}

\begin{proof}
Posons $U ^\dag _\alpha := P ^\dag _\alpha \setminus T _\alpha$,
$U ^\dag _{\alpha \beta}:= P ^\dag _\alpha \times P ^\dag _\beta  \setminus T _{\alpha \beta}$
et fixons
 $Y ^\dag _\alpha \hookrightarrow U ^\dag _\alpha$,
 $Y ^\dag _{\alpha \beta} \hookrightarrow U ^\dag _{\alpha \beta}$
 des relèvements respectifs de
 $Y _\alpha \subset U _\alpha$, $Y _{\alpha \beta }\subset U _{\alpha \beta}$.
Soient
$((E _\alpha) _{\alpha \in \Lambda}, (\eta _{  \alpha \beta}) _{\alpha,\beta  \in \Lambda})
\in
F\text{-}\mathrm{Isoc} ^{\dag}  (Y,\, (Y _\alpha, P ^\dag _\alpha, T_\alpha,  X_\alpha) _{\alpha \in \Lambda}/K).$
Pour tous $\alpha,\beta \in \Lambda$,
notons $\E _\alpha :=\sp _{Y ^\dag _\alpha \hookrightarrow U ^\dag _\alpha, T _\alpha,+} (E _\alpha)
\in F\text{-}\mathrm{Isoc} ^{\dag \dag} (Y _\alpha /K)$ (voir \ref{dag=dagdag-affine})
et $\theta  _{\alpha \beta}$
l'isomorphisme
$j _2 ^{\alpha \beta *} (\E _\beta)
  \riso j _1 ^{\alpha \beta *} (\E _\alpha)$
défini via le diagramme commutatif :
\small
\begin{equation}
\label{recollcompsp+-diag}
\xymatrix @R =0,3cm {
{\R \underline{\Gamma} ^\dag _{X _{{\alpha \beta}}}
\circ (\hdag T _{{\alpha \beta}})\circ p _2 ^{\alpha \beta !} (\sp _{Y ^\dag _\beta \hookrightarrow U ^\dag _\beta, T _\beta,+} (E _\beta))}
\ar@{.>}[r] ^-\sim _-{\theta _{  \alpha \beta}}
\ar[d] ^-\sim
&
{\R \underline{\Gamma} ^\dag _{X _{{\alpha \beta}}}
\circ (\hdag T _{{\alpha \beta}})\circ p _1 ^{\alpha \beta !}
(\sp _{Y ^\dag _\alpha \hookrightarrow U ^\dag _\alpha, T _\alpha,+} (E _\alpha))}
\ar[d] ^-\sim
\\
{\sp _{Y ^\dag _{\alpha \beta} \hookrightarrow U ^\dag _{\alpha \beta}, T _{\alpha \beta},+} j _2 ^{\alpha \beta *}
(E _\beta)}
\ar[r] ^-\sim
_-{\sp _{Y ^\dag _{\alpha \beta} \hookrightarrow U ^\dag _{\alpha \beta}, T _{\alpha \beta},+} (\eta _{  \alpha \beta})}
&
{\sp _{Y ^\dag _{\alpha \beta} \hookrightarrow U ^\dag _{\alpha \beta}, T _{\alpha \beta},+} j _1 ^{\alpha \beta *} (E _\alpha),}
}
\end{equation}
\normalsize
où les isomorphismes verticaux résultent de \cite[5.2.3 et 5.2.6]{caro_devissge_surcoh}
(par fidélité, ces isomorphismes commutent à Frobenius car c'est le cas en dehors de $T _{\alpha \beta}$).
Par transitivité des isomorphismes verticaux de \ref{recollcompsp+-diag}, on vérifie
la condition de cocycle de la famille $(\theta _{  \alpha \beta}) _{\alpha,\beta  \in \Lambda}$.
D'où
le foncteur
$\sp _+$
:
$F\text{-}\mathrm{Isoc} ^{\dag}  (Y,\, (Y _\alpha, P ^\dag _\alpha, T_\alpha,  X_\alpha) _{\alpha \in \Lambda}/K)
\rightarrow
F\text{-}\mathrm{Isoc} ^{\dag \dag}  (Y,\, (Y _\alpha, \PP _\alpha, T_\alpha,  X_\alpha) _{\alpha \in \Lambda}/K)$
défini par
$((E _\alpha) _{\alpha \in \Lambda}, (\eta _{  \alpha \beta}) _{\alpha,\beta  \in \Lambda})
\mapsto
((\E _\alpha) _{\alpha \in \Lambda}, (\theta _{  \alpha \beta}) _{\alpha,\beta  \in \Lambda}) $.
Il résulte de \ref{eqcataffine} que $\sp _{(Y _\alpha, P ^\dag _\alpha , T_\alpha,  X _\alpha) _{\alpha \in \Lambda},+}$
est une équivalence de catégories.
Or, via \ref{ind-isocsurcvgen-vide}
et \ref{dag=dagafflisse}, on obtient
$F\text{-}\mathrm{Isoc} ^{\dag}  (Y/K) \cong
F\text{-}\mathrm{Isoc} ^{\dag}  (Y,\, (Y _\alpha, P ^\dag _\alpha, T_\alpha,  X_\alpha) _{\alpha \in \Lambda}/K)$.
D'où le résultat par composition de ces équivalences.
\hfill \hfill \qed \end{proof}

\section{Surcohérence différentielle des $F$-isocristaux surconvergents sur les variétés lisses}

\subsection{Contagiosité de la surcohérence}

Précisons et étendons d'abord la remarque \cite[2.2.7]{caro_surcoherent} :
\begin{prop}
\label{Rsp*comm}
Soient $\PP$ un $\V$-schéma formel lisse, $T$ un diviseur de $P$, $\U := \PP \setminus T$.
  Soit $E$ un isocristal sur $U$ surconvergent le long de $T$, i.e., un $j ^\dag _T \O _{\PP _K}$-module cohérent muni
  d'une connexion surconvergente.
  Pour tous sous-schémas fermés $Z$ et $Z'$ de $P$,
  on dispose alors d'un isomorphisme canonique :
\begin{equation}
\label{Rsp*comm-iso}
  \R \sp _* \underline{\Gamma} ^\dag _{]Z '[ } j ^\dag _Z (E) \riso
  \R \underline{\Gamma} ^\dag _{Z'} (\hdag Z) \sp _* (E).
\end{equation}
\end{prop}

\begin{proof}
\'Etablissons d'abord l'existence d'un isomorphisme (non nécessairement canonique) de la forme \ref{Rsp*comm-iso}.
Via les triangles distingués
\begin{gather}
\notag
  \R \sp _*  \underline{\Gamma} ^\dag _{]Z '[ } ( \underline{\Gamma} ^\dag _{]Z [ }(E) )
\rightarrow
\R \sp _*  \underline{\Gamma} ^\dag _{]Z '[ } (E)
\rightarrow
\R \sp _*  \underline{\Gamma} ^\dag _{]Z '[ } (j ^\dag _Z (E) )
\rightarrow +1,
\\
\notag
\R \underline{\Gamma} ^\dag _{Z'} (\R \underline{\Gamma} ^\dag _{Z} \sp _* (E))
\rightarrow
\R \underline{\Gamma} ^\dag _{Z'} (\sp _* (E))
\rightarrow
\R \underline{\Gamma} ^\dag _{Z'} ((\hdag Z) \sp _* (E))
\rightarrow
+1,
\end{gather}
on se ramène au cas où $Z$ est vide. Prouvons alors le lemme par récurrence
sur le nombre minimum de diviseurs d'intersection $Z'$. Lorsque $Z' $ est un diviseur,
$\R \sp _*  j ^\dag _{Z'} (E) \liso \sp _*  j ^\dag _{Z'} (E) \riso (\hdag Z') \sp _*   (E)$.
En utilisant le triangle de localisation en $Z'$ de $\sp _* (E)$ (voir \cite[2.2.6]{caro_surcoherent})
et en appliquant $\R \sp _*  $ à la suite exacte
$0 \rightarrow \underline{\Gamma} ^\dag _{]Z '[ } (E) \rightarrow
E \rightarrow j ^\dag _{Z'} (E) \rightarrow 0$ (voir \cite[2.1.6]{Berig}),
il en résulte
$\R \underline{\Gamma} ^\dag _{Z'} (\sp _* (E)) \riso \R \sp _*  \underline{\Gamma} ^\dag _{]Z '[ } (E)$.
De plus, si $T _1, \dots, T _r$ sont des diviseurs d'intersection $Z'$ avec $r \geq 2$,
en appliquant $\R \sp _*$ à une suite exacte
de Mayer Vietoris (cela découle de \cite[2.1.7]{Berig}), on obtient le triangle distingué :
\small
\begin{equation}
  \notag
\R \sp _*  \underline{\Gamma} ^\dag _{]Z '[ } (E)
\rightarrow
\R \sp _*  \underline{\Gamma} ^\dag _{]T _1 \cap \dots \cap T _{r-1}[ } (E)
\oplus
\R \sp _*  \underline{\Gamma} ^\dag _{]T _r[ } (E)
\rightarrow
\R \sp _*  \underline{\Gamma} ^\dag _{](T _1 \cap \dots \cap T _{r-1} )\cup T _r[ } (E)
\rightarrow
+1.
\end{equation}
\normalsize
Comme on bénéficie du triangle distingué de Mayer Vietoris (voir \cite[2.2.16]{caro_surcoherent})
\small
\begin{equation}
  \notag
\R \underline{\Gamma} ^\dag _{Z'} \sp _* (E)
\rightarrow
\R \underline{\Gamma} ^\dag _{T _1 \cap \dots \cap T _{r-1}} \sp _* (E)
\oplus
\R \underline{\Gamma} ^\dag _{T _r} \sp _* (E)
\rightarrow
\R \underline{\Gamma} ^\dag _{(T _1 \cap \dots \cap T _{r-1} )\cup T _r} \sp _* (E)
\rightarrow
+1,
\end{equation}
\normalsize
on conclut la récurrence.

Prouvons à présent que l'existence d'un isomorphisme de la forme \ref{Rsp*comm-iso} nous permet
d'en construire un canoniquement.
En appliquant le foncteur $(\hdag Z) \circ \R \sp _*$ à la suite exacte
$0 \rightarrow \underline{\Gamma} ^\dag _{]Z [ } (E) \rightarrow
E \rightarrow j ^\dag _Z (E) \rightarrow 0$,
on obtient un triangle distingué dont le terme de gauche est
nul grâce à (l'isomorphisme de la forme) \ref{Rsp*comm-iso}.
Il en dérive l'isomorphisme canonique
$(\hdag Z) \circ \R \sp _* (E) \riso (\hdag Z) \circ \R \sp _* (j ^\dag _Z (E) ) $.
Par \ref{Rsp*comm-iso} (et \cite[2.2.14]{caro_surcoherent}), le morphisme canonique
$\R \sp _* (j ^\dag _Z (E) )
\rightarrow
(\hdag Z) \circ \R \sp _* (j ^\dag _Z (E) )$
est un isomorphisme. D'où l'isomorphisme canonique
$\R \sp _* (j ^\dag _Z (E) )  \riso (\hdag Z) \circ \R \sp _* (E)$
puis
$\R \underline{\Gamma} ^\dag _{Z'}  \R \sp _* (j ^\dag _Z (E) )
\riso
\R \underline{\Gamma} ^\dag _{Z'} (\hdag Z) \circ \R \sp _* (E)$.
De même, en appliquant le foncteur
$\R \underline{\Gamma} ^\dag _{Z'} \circ \R \sp _* $
à la suite exacte
$0 \rightarrow \underline{\Gamma} ^\dag _{]Z ' [ } j ^\dag _Z (E)  \rightarrow
j ^\dag _Z (E)  \rightarrow j ^\dag _{Z \cup Z'} (E) \rightarrow 0$,
on obtient les isomorphismes canoniques :
$\R \sp _* (\underline{\Gamma} ^\dag _{]Z ' [ } j ^\dag _Z (E)) \liso
\R \underline{\Gamma} ^\dag _{Z'} \circ  \R \sp _* (\underline{\Gamma} ^\dag _{]Z ' [ } j ^\dag _Z (E)) \riso
\R \underline{\Gamma} ^\dag _{Z'} \circ \R \sp _* (j ^\dag _Z (E)  ). $
Il en résulte par composition la construction de l'isomorphisme canonique \ref{Rsp*comm-iso}.
\hfill \hfill \qed \end{proof}

\begin{prop}
\label{lemm-isoc-surcoh}
Soient $\PP$ un $\V$-schéma formel lisse,
$\E\in \smash{\underset{^{\longrightarrow}}{LD}} ^{\mathrm{b}} _{\Q,\mathrm{qc}}
( \smash{\widehat{\D}} _{\PP} ^{(\bullet)})$, $T$ un diviseur de $P$,
$T _1,\dots ,T _r$ des diviseurs de $P$ d'intersection $T$.

$(a)$ Pour tout $\alpha = 1,\dots ,r$,
supposons qu'il existe un ($F$-)isocristal $E _\alpha$ sur $P \setminus T _\alpha$ surconvergent le long de $T _\alpha $
et un isomorphisme ($F$-)$\D ^\dag _{\PP} (\hdag T) _\Q$-linéaire
$\E (\hdag T _\alpha) \riso \sp _*  (E _\alpha)$ (on omet d'indiquer le foncteur $\underset{\longrightarrow}{\lim}$).
Il existe alors un ($F$-)isocristal $E $ sur $P \setminus T $ surconvergent le long de $T  $
tel que $j ^\dag _{T _\alpha} E \riso E _\alpha$
qui induit une famille compatible d'isomorphismes ($F$-)$\D ^\dag _{\PP} (\hdag T) _\Q$-linéaires
$\E (\hdag T _\alpha) \riso \sp _*  ( j ^\dag _{T_ \alpha} E)$ (i.e., ils induisent canoniquement les mêmes
isomorphismes $\E (\hdag T _{\alpha \cup \beta }) \riso \sp _*  ( j ^\dag _{T_ \alpha \cup T_ \beta} E)$).

$(b)$ Soit $E $ un ($F$-)isocristal sur $P \setminus T $ surconvergent le long de $T  $ tel que,
pour tout $\alpha = 1,\dots ,r$, il existe un isomorphisme ($F$-)$\D ^\dag _{\PP} (\hdag T _\alpha) _\Q$-linéaire :
$\E (\hdag T _\alpha) \riso \sp _*  ( j ^\dag _{T_ \alpha} E)$.
Il existe alors un isomorphisme ($F$-)$\D ^\dag _{\PP} (\hdag T) _\Q$-linéaire
de la forme
$\E (\hdag T) \riso \sp _* (E)$.
Si les isomorphismes
$\E (\hdag T _\alpha) \riso \sp _*  ( j ^\dag _{T_ \alpha} E)$ sont en outre compatibles,
il existe alors un unique isomorphisme ($F$-)$\D ^\dag _{\PP} (\hdag T) _\Q$-linéaire
$\E (\hdag T) \riso \sp _* (E)$ induisant canoniquement pour tout $\alpha = 1,\dots ,r$ les isomorphismes
$\E (\hdag T _\alpha) \riso \sp _*  ( j ^\dag _{T_ \alpha} E)$.
\end{prop}

\begin{proof}
Prouvons d'abord $(a)$.
On construit les isomorphismes canoniques
$j ^\dag _{T _\alpha} E _\beta \riso j ^\dag _{T _\beta} E _\alpha$
via le diagramme commutatif suivant :
\begin{equation}
  \label{lemm-isoc-surcoh-diag0}
  \xymatrix @R=0,3cm {
{(\hdag T _\beta)  (\hdag T _\alpha)  (\E )}
\ar[r] _-\sim
&
{(\hdag T _\beta)  ( \sp _*  (E _\alpha) )}
\ar[r] _-\sim
&
{ \sp _*  (j ^\dag _{T _\beta} E _\alpha) }
\\
{  (\hdag T _\alpha) (\hdag T _\beta) (\E )}
\ar[r] _-\sim \ar[u] _-\sim
&
{(\hdag T _\alpha)  ( \sp _*  (E _\beta) )}
\ar[r] _-\sim
&
{ \sp _*  (j ^\dag _{T _\alpha} E _\beta) .}
\ar@{.>}[u] _-\sim
}
\end{equation}
Comme ces isomorphismes vérifient la condition de
cocycle,
il découle alors de \cite[2.1.12]{Berig} l'existence
d'un ($F$-)isocristal $E $ sur $P \setminus T $ surconvergent le long de $T  $
tel que $E _\alpha \riso j ^\dag _{T _\alpha} E$.
La compatibilité est tautologique (on complète le diagramme \ref{lemm-isoc-surcoh-diag0} en ajoutant
$\E (\hdag T _{\alpha \cup \beta })$ et $ \sp _*  ( j ^\dag _{T_ \alpha \cup T_ \beta} E)$).

Traitons à présent $(b)$. Pour cela,
vérifions par récurrence sur l'entier $s$ tel que $1\leq s \leq r$ l'assertion suivante :
{\og pour tous diviseurs $T '_1,\dots, T' _s$ de $P$ contenant respectivement $T _1,\dots, T_s$,
$\E (\hdag T' _1 \cap \dots \cap T '_s) \riso \R \sp _* (j ^\dag _{T' _1 \cap \dots \cap T '_s} E)$\fg}.

Lorsque $s =1$,
  $\E (\hdag T _1) \riso \sp _* j ^\dag _{T _1} (E)$ implique
      $\E (\hdag T '_1) \riso (\hdag T '_1) \circ \sp _* j ^\dag _{T _1} (E)\newline
      \underset{\tiny{\ref{Rsp*comm-iso}}}{\riso} \sp _* j ^\dag _{T '_1} (E)$.
    Supposons à présent $s \geq 2$. En notant, pour $\beta =1, \dots , s -1$, $T '' _\beta := T ' _\beta \cup T ' _s$,
    on dispose des triangles distingués de Mayer Vietoris (voir \cite[2.2.16]{caro_surcoherent} pour le premier) :
\small
    \begin{gather}
    \notag
      \E (\hdag T ' _1 \cap \dots \cap T ' _s)
      \rightarrow
      \E (\hdag T ' _1 \cap \dots \cap T ' _{s-1}) \oplus \E (\hdag T ' _s)
      \rightarrow
      \E (\hdag T '' _1 \cap \dots \cap T '' _{s-1})
      \rightarrow
      +1,\\ \notag
      \R \sp _* j ^\dag _{T ' _1 \cap \dots \cap T ' _s}(E)
      \rightarrow
      \R \sp _* j ^\dag _{T ' _1 \cap \dots \cap T ' _{s-1}}(E) \oplus \R \sp _* j ^\dag _{T ' _{s}}(E)
      \rightarrow
      \R \sp _* j ^\dag _{T '' _1 \cap \dots \cap T '' _{s-1}}(E)
      \rightarrow
      +1.
    \end{gather}
    \normalsize
    On conclut ainsi la récurrence. En particulier, le cas $r =s$ fournit un isomorphisme
    $\E (\hdag T) \riso \sp _* (E)$.

Supposons à présent la famille d'isomorphismes $\E (\hdag T _\alpha) \riso \sp _*  ( j ^\dag _{T_ \alpha} E)$
compatible. Avec \cite[6.1.4]{caro_devissge_surcoh}
(l'hypothèse séparé est ici inutile car elle ne sert qu'à construire par recollement les foncteurs de la forme
$\sp _{X \hookrightarrow \PP,T,+}$),
il en découle que
    $\E (\hdag T) , \sp _* (E) \in (F\text{-})\mathrm{Isoc} ^{\dag \dag} ( \PP ,\, T, \,P/K )$.
    Comme le foncteur $(\hdag T _\alpha)$ :
    $(F\text{-})\mathrm{Isoc} ^{\dag \dag} ( \PP ,\, T, \,P/K ) \rightarrow
    (F\text{-})\mathrm{Isoc} ^{\dag \dag} ( \PP ,\, T_\alpha, \,P/K )$ est pleinement fidèle
    (voir \ref{plfdttssfrob}), il existe un unique isomorphisme $\E (\hdag T) \riso \sp _* (E)$
    induisant $\E (\hdag T _\alpha) \riso \sp _*  ( j ^\dag _{T_ \alpha} E)$
    via les isomorphismes canoniques de la forme $\E (\hdag T) (\hdag T _\alpha) \riso (\hdag T _\alpha)$
    et $\sp _* (E) (\hdag T _\alpha) \riso \sp _* ( j^\dag _{T _\alpha} E)$.
    La compatibilité entraîne que cet isomorphisme $\E (\hdag T) \riso \sp _* (E)$
    ne dépend pas du choix de $\alpha$.
\hfill \hfill \qed \end{proof}

\begin{nota}
\label{nota-isoc-surcoh}
  Soient $\PP$ un $\V$-schéma formel propre et lisse,
  $u$ : $X \hookrightarrow P$ une immersion fermée, $T$ un diviseur de $P$,
  $\U$ l'ouvert de $\PP$ complémentaire de $T$.
  On suppose $Y := X \setminus T$ lisse sur $k$.
\end{nota}

\begin{theo}
\label{isocsurcohgen-prel}
Avec les notations de \ref{nota-isoc-surcoh},
soient $\E$ un $F$-$\D ^\dag _{\PP} (\hdag T) _{\Q}$-module cohérent à support dans $X$,
  $T _1, \dots , T _r$ des diviseurs de $P$ d'intersection $T$.

Si, pour tout $\alpha =1,\dots , r$,
$\E (\hdag T _\alpha) \in F\text{-}\mathrm{Isoc} ^{\dag \dag} (\PP , T _\alpha , X /K)$ alors
$\E \in F\text{-}\mathrm{Isoc} ^{\dag \dag} (\PP , T , X /K)$.
\end{theo}

\begin{proof}
  Il ne coûte rien de supposer $Y$ dense dans $X$ et $X$ irréductible (on peut utiliser \cite[6.3.3]{caro_devissge_surcoh}).
  D'après \cite[6.3.1]{caro_devissge_surcoh},
  on dispose d'un diagramme de la forme \ref{diag631dejong}
  (on conserve les notations et hypothèses de \ref{diag631dejong}).
Notons $\E ':=\R \underline{\Gamma} ^\dag _{X'} f ^! (\E )$,
$T ' _\alpha:= f ^{-1} (T _\alpha) $, $Y '_\alpha := Y '\setminus T '_\alpha$.
Le fait que $T ' \cap X'$ soit un diviseur de $X'$, $T ' \subset T '_\alpha$ et $X'$ soit lisse implique que
les $T ' _\alpha \cap X'$ sont aussi des diviseurs de $X'$.
Notons alors
$E ' _\alpha$ le $F$-isocristal surconvergent sur $Y ' _\alpha$ associé à $\E '(\hdag T ' _\alpha)$
(via $\sp _{X' \hookrightarrow \PP', T' _\alpha,+}$).
Comme $X'$ est lisse,
le foncteur $\sp _{ X' \hookrightarrow \PP', T', +}$
induit l'équivalence de catégories :
$F\text{-}\mathrm{Isoc} ^{\dag} (Y '/K) \cong
F\text{-}\mathrm{Isoc} ^{\dag \dag} (\PP' , T ', X '/K)$.
Via la description de l'image essentielle de
$\sp _{ X' \hookrightarrow \PP', T', +}$ (voir \cite[4.1.9]{caro-construction}),
le fait que $\E ' \in F\text{-}\mathrm{Isoc} ^{\dag \dag} (\PP' , T ', X '/K)$
est donc local en $\PP'$. Pour le vérifier, on peut ainsi supposer qu'il existe un
relèvement $u'$ : $\X '\hookrightarrow \PP'$ de $X' \hookrightarrow \PP'$.
Dans ce cas,
$(\hdag T' _\alpha \cap X') (u ^{\prime !} (\E') )\riso
u ^{\prime !}   (\E'(\hdag T' _\alpha ) )
\riso
\sp _* (E' _\alpha)$, où
$\sp$ : $\X _K \rightarrow \X$ est le morphisme de spécialisation.
Par \ref{lemm-isoc-surcoh},
il existe alors un unique
$F$-isocristal $E'$ surconvergent sur $Y'$ tel que
$u ^{\prime !} (\E ')\riso \sp _* (E')$. Il en résulte que
$\E' \in F\text{-}\mathrm{Isoc} ^{\dag \dag} (\PP' , T ', X '/K)$.
De la même façon, on établit via \ref{lemm-isoc-surcoh} que
$\E |\U$ est dans l'image essentielle de $\sp _{Y \hookrightarrow \U, +}$.
Il en découle que le $F$-complexe $\DD _T (\E)$ est isomorphe à
un $F$-$\D ^\dag _{\PP} (\hdag T) _{\Q}$-module cohérent à support dans $X$
(car c'est le cas en dehors de $T$ puis on utilise \cite[4.3.12]{Be1}).
Comme, pour tout $1 \leq \alpha \leq r$,
$(\hdag T _\alpha) (\DD _T \E ) \in F\text{-}\mathrm{Isoc} ^{\dag \dag} (\PP , T _\alpha , X /K)$,
on obtient de même que
$\R \underline{\Gamma} ^\dag _{X'} f_T ^! (\DD _T \E )
\in
F\text{-}\mathrm{Isoc} ^{\dag \dag} (\PP' , T ', X '/K)$.
En outre,
on vérifie l'isomorphisme :
$\DD _{T'} \R \underline{\Gamma} ^\dag _{X'} f ^! (\DD _T \E ) \riso
\R \underline{\Gamma} ^\dag _{X'} f ^! (\E )  $
(comme $Y$ est lisse, cela est aisé en dehors $T'$ puis on utilise \cite[6.5.2]{caro_devissge_surcoh}).
Or, on dispose des morphismes
\begin{gather}
\notag
f _{T,+} \R \underline{\Gamma} ^\dag _{X'} f_T ^! (\E ) \rightarrow \E, \\
\notag  \E \rightarrow f _{T,+} \DD _{T'} \R \underline{\Gamma} ^\dag _{X'} f_T ^! (\DD _T \E ),
\end{gather}
le premier se construisant par adjonction, le second s'en déduisant via les isomorphismes de bidualité
(voir \cite[1.3.6]{virrion})
et de dualité relative
(voir \cite[1.2.7]{caro_courbe-nouveau}).
D'où les flèches :
\begin{equation}
  \label{comp-fact-dir}
\E \rightarrow f _{T,+} \R \underline{\Gamma} ^\dag _{X'} f_T ^! (\E ) \rightarrow \E.
\end{equation}
Comme notamment $\E |\U$ est dans l'image essentielle de
$\sp _{Y \hookrightarrow \U, +}$, on vérifie que
le composé de \ref{comp-fact-dir} est un isomorphisme en dehors de $T$.
Par \cite[4.3.12]{Be1}, il en résulte que celui-ci est un isomorphisme et donc que
$\E$ est un facteur directe de $f _{T,+} \R \underline{\Gamma} ^\dag _{X'} f_T ^! (\E )$.
On conclut via \cite[6.1.3]{caro_devissge_surcoh}.
\hfill \hfill \qed \end{proof}

\begin{prop}
  \label{lemm-eqcat-rel}
Avec les notations de \ref{nota-isoc-surcoh},
  soit $T _1$ un diviseur de $P$ contenant $T$ tel que
  $Y _1:= X \setminus T _1$ soit affine.
  Soient $E \in F\text{-}\mathrm{Isoc} ^{\dag } ( Y/K)$
  et
  $\E _1 \in F\text{-}\mathrm{Isoc} ^{\dag \dag} ( \PP ,\, T_1, \,X/K )$ tel que
  $\sp _{Y _1,+}( j^\dag _{T _1} E) \riso \E _1$ (notation de \ref{eqcataffine}).
  Alors, $\E _1 $ et $\DD _{T } (\E _1)$ sont $\D ^\dag _{\PP} (\hdag T) _{\Q}$-surcohérents.
\end{prop}

\begin{proof}
  On peut supposer $X$ irréductible et $T \not \supset X$.
  Avec ses notations et hypothèses,
  on dispose alors du diagramme \ref{diag631dejong}.
Notons $\E ' _1:=\R \underline{\Gamma} ^\dag _{X'} f ^! (\E _1 )$,
$T ' _1:= f ^{-1} (T _1) $, $Y '_1 := Y '\setminus T '_1$,
$E ' $ le $F$-isocristal surconvergent sur $Y '$ égal à $b ^* (E)$.
L'isomorphisme $\sp _{Y _1,+}( j^\dag _{T _1} E) \riso \E _1$
induit le suivant
$\sp _{X ' \hookrightarrow \PP', T' _1,+}( j^\dag _{T ' _1} E')  \riso \E ' _1$ (en effet,
grâce à la  pleine fidélité de \cite[6.5.2]{caro_devissge_surcoh},
il suffit de le valider en dehors de $T'_1$, ce qui est immédiat car tout est lisse en dehors des diviseurs).
Il en résulte que
$(\hdag T ' _1 )\circ \sp _{X ' \hookrightarrow \PP', T' ,+}( E')  \riso \E ' _1$.
Par \cite[6.1.4]{caro_devissge_surcoh},
il en dérive que
$\E ' _1$ et $\DD _{T' } (\E ' _1)$ sont
$\D ^\dag _{\PP'} (\hdag T') _{\Q}$-surcohérents.
Or, d'après \cite[6.3.1]{caro_devissge_surcoh},
$\E  _1$
est un facteur direct de $f _+ (\E ' _1)$.
Il en résulte, via le théorème de dualité relative
(voir \cite[1.2.7]{caro_courbe-nouveau}),
que $\DD _{T} (\E _1)$ est
un facteur direct de $f _+ \circ \DD _{T' } (\E ' _1)$.
On conclut avec la stabilité de la surcohérence par l'image directe d'un morphisme propre (voir \cite[3.1.9]{caro_surcoherent}).
\hfill \hfill \qed \end{proof}

\begin{rema}
\label{plfdttalpha}
  Avec les notations de \ref{nota-isoc-surcoh},
  soient $T _1$ un diviseur de $P$ contenant $T$,
  $Y _1:= X \setminus T _1$.
  Le foncteur $(\hdag T_1)$ :
  $F\text{-}\mathrm{Isoc} ^{\dag \dag} ( \PP ,\, T, \,X/K )
  \rightarrow
  F\text{-}\mathrm{Isoc} ^{\dag \dag} ( \PP ,\, T_1, \,X/K )$ est pleinement fidèle.
  En effet, on dispose du diagramme commutatif (à isomorphisme près) :
  $$ \xymatrix @R=0,3cm  {
  {F\text{-}\mathrm{Isoc} ^{\dag \dag} ( \PP ,\, T, \,X/K )}
  \ar[r]
  \ar[d] _-{\rho _{Y }}
  &
  {F\text{-}\mathrm{Isoc} ^{\dag \dag} ( \PP ,\, T_1, \,X/K )}
  \ar[d] _-{\rho _{Y _1}}
  \\
  {F\text{-}\mathrm{Isoc} ^{\dag } ( Y/K)}
  \ar[r]
  &
  {F\text{-}\mathrm{Isoc} ^{\dag } ( Y _1/K),}
  }$$
dont les flèches verticales et du bas sont pleinement fidèles
(voir respectivement
\cite[6.5.3]{caro_devissge_surcoh} et \cite[4.1.2]{tsumono}).
\end{rema}

\subsection{$F$-isocristaux surcohérents sur les variétés lisses}

\begin{nota}
  \label{notaprop1etape}
  Soient $\PP$ un $\V$-schéma formel propre et lisse,
  $u$ : $X \hookrightarrow P$ une immersion fermée, $T$ un diviseur de $P$,
  $\U$ l'ouvert de $\PP$ complémentaire de $T$.
  On suppose $P$ intègre, $Y := X \setminus T$ intègre et lisse sur $k$.
De plus, soient $(\PP _\alpha) _{\alpha \in \Lambda}$ un recouvrement fini d'ouverts affines de $\PP$,
$Y _{\alpha }:=P _\alpha \cap Y$ et
$T _\alpha := T \cup (P \setminus P _\alpha)$ (ce dernier est un diviseur de $P$ vérifiant
$X \setminus T _\alpha = Y _\alpha$).

\end{nota}

\begin{rema}
\label{remaprop1etape}
  Avec les notations \ref{notaprop1etape}, soit
  $((\E _\alpha) _{\alpha \in \Lambda}, (\theta _{  \alpha \beta}) _{\alpha,\beta  \in \Lambda}) \in
F\text{-}\mathrm{Isoc} ^{\dag \dag} ( Y,\, (Y _\alpha, \PP , T_\alpha,  X) _{\alpha \in \Lambda}/K)$.
Notons respectivement $p_1$ et $p_2$ les projections à gauche et à droite
$\PP \times \PP \rightarrow \PP$,
$T _{\alpha \beta} := p _1 ^{-1} (T _\alpha) \cup p _2 ^{-1} (T _\beta)$
et $\delta$ : $\PP \hookrightarrow \PP \times \PP$ l'immersion fermée diagonale.
En appliquant
$\delta ^!$ à
$\theta _{\alpha \beta}$ :
$\R \underline{\Gamma} ^\dag _{X }
\circ (\hdag T _{\alpha \beta})\circ p _2 ^{ !}
(\E _\beta)
\riso
\R \underline{\Gamma} ^\dag _{X }
\circ (\hdag T _{\alpha \beta})\circ p _1 ^{ !}
(\E _\alpha)$, on obtient
l'isomorphisme
$\vartheta _{\alpha \beta}$ :
$(\hdag T _\alpha)(\E _\beta)
\riso
(\hdag T _\beta )  ( \E _\alpha)$. Ainsi,
$F\text{-}\mathrm{Isoc} ^{\dag \dag} ( Y,\, (Y _\alpha, \PP , T_\alpha,  X) _{\alpha \in \Lambda}/K)$
est équivalente à la catégorie dont les objets sont de la forme
$((\E _\alpha) _{\alpha \in \Lambda}, (\vartheta _{  \alpha \beta}) _{\alpha,\beta  \in \Lambda}) $,
où
$\E _\alpha \in F\text{-}\mathrm{Isoc} ^{\dag \dag} ( \PP ,\, T_\alpha, \,X/K )$
et $\vartheta _{\alpha \beta } $ :
$(\hdag T _\alpha)(\E _\beta)
\riso
(\hdag T _\beta )  ( \E _\alpha)$
({\it la donnée de recollement})
est un isomorphisme vérifiant la condition de cocycle
$(\hdag T _\beta) (\vartheta _{\alpha \gamma })=
(\hdag T _\gamma) (\vartheta_{  \alpha \beta} )
\circ (\hdag T _\alpha) ( \vartheta _{ \beta \gamma })$, dont
les morphismes de la forme
$((\E '_\alpha) _{\alpha \in \Lambda}, (\vartheta '_{  \alpha \beta}) _{\alpha,\beta  \in \Lambda})  \rightarrow
((\E _\alpha) _{\alpha \in \Lambda}, (\vartheta _{  \alpha \beta}) _{\alpha,\beta  \in \Lambda})$
sont les morphismes
$\E ' _\alpha \rightarrow \E _\alpha$ de
$F\text{-}\mathrm{Isoc} ^{\dag \dag} ( \PP ,\, T_\alpha, \,X/K )$
commutant aux données de recollement.

\end{rema}

\begin{prop}
  \label{eqcat-rel}
Avec les notations de \ref{notaprop1etape},
on dispose de l'équivalence canonique de catégories  :
$$F\text{-}\mathrm{Isoc} ^{\dag \dag} (\PP, \, T , \, X/K)
\cong
F\text{-}\mathrm{Isoc} ^{\dag \dag} ( Y,\, (Y _\alpha, \PP , T_\alpha,  X) _{\alpha \in \Lambda}/K).$$
\end{prop}

\begin{proof}
On construit le foncteur
$\mathcal{L}oc$ :
$F\text{-}\mathrm{Isoc} ^{\dag \dag} (\PP, \, T , \, X/K)\rightarrow \linebreak[4]
F\text{-}\mathrm{Isoc} ^{\dag \dag} ( Y,\, (Y _\alpha, \PP , T_\alpha,  X) _{\alpha \in \Lambda}/K)$
de la manière suivante :
on associe à $\E \in F\text{-}\mathrm{Isoc} ^{\dag \dag} (\PP, \, T , \, X/K)$,
$( \E (\hdag T _\alpha), \vartheta _{\alpha \beta})\in
F\text{-}\mathrm{Isoc} ^{\dag \dag} ( Y,\, (Y _\alpha, \PP , T_\alpha,  X) _{\alpha \in \Lambda}/K)$, où
$\vartheta_{\alpha \beta}$ est l'isomorphisme canonique :
$(\hdag T _\alpha)(\hdag T _\beta) (\E)
\riso
(\hdag T _\beta) (\hdag T _\alpha) (\E )$ (voir la remarque \ref{remaprop1etape}).
Il découle de \ref{plfdttalpha} que $\mathcal{L}oc$ est pleinement fidèle.
Prouvons à présent via les étapes $(a)$ et $(b)$ ci-dessous son essentielle surjectivité.

$(a)$ Avec ses notations et hypothèses, on bénéficie du diagramme \ref{diag631dejong}.
On pose aussi $\U _\alpha := \PP \setminus T _\alpha$,
 $T ' _\alpha:= f ^{-1} (T _\alpha) $,
 $\U '_\alpha := \PP '\setminus T '_\alpha$, $Y '_\alpha := Y '\setminus T '_\alpha$.
On dispose alors du foncteur canonique
$F\text{-}\mathrm{Isoc} ^{\dag \dag} ( Y,\, (Y _\alpha, \PP , T_\alpha,  X) _{\alpha \in \Lambda}/K)
\rightarrow
F\text{-}\mathrm{Isoc} ^{\dag \dag} ( Y',\, (Y '_\alpha, \PP ', T'_\alpha,  X') _{\alpha \in \Lambda}/K)$
via
$(\E _\alpha, \theta _{\alpha \beta})
\mapsto
(\R \underline{\Gamma} ^\dag _{X'} f_T ^! ( \E _\alpha), \theta '_{\alpha \beta})$,
où $\theta '_{\alpha \beta}$ est l'isomorphisme canonique induit par
$\theta _{\alpha \beta}$ (on utilise les isomorphismes canoniques \cite[2.2.18]{caro_surcoherent}).
On bénéficie aussi du foncteur canonique
$F\text{-}\mathrm{Isoc} ^{\dag \dag} ( Y,\, (Y _\alpha, \PP , T_\alpha,  X) _{\alpha \in \Lambda}/K)
\rightarrow
F\text{-}\mathrm{Isoc} ^{\dag \dag} ( Y,\, (Y _\alpha, \U _\alpha) _{\alpha \in \Lambda}/K)$
via
$(\E _\alpha, \theta _{\alpha \beta})
\mapsto
(\E _\alpha |\U _\alpha , \theta _{\alpha \beta} |\U _\alpha \times \U _\beta )$, de même
avec des primes.
Par \ref{congcascompliss}, on a
$F\text{-}\mathrm{Isoc} ^{\dag} ( Y'/K)
\cong
F\text{-}\mathrm{Isoc} ^{\dag \dag} ( Y',\, (Y '_\alpha, \PP ', T'_\alpha,  X') _{\alpha \in \Lambda}/K)$,
$F\text{-}\mathrm{Isoc}  ( Y/K)
\cong \linebreak[3]
F\text{-}\mathrm{Isoc} ^{\dag \dag} ( Y,\, (Y _\alpha, \U _\alpha) _{\alpha \in \Lambda}/K)$,
$F\text{-}\mathrm{Isoc}  ( Y'/K)
\cong
F\text{-}\mathrm{Isoc} ^{\dag \dag} ( Y',\, (Y' _\alpha, \U '_\alpha) _{\alpha \in \Lambda}/K)$.
D'où le foncteur :
\small
\begin{equation}
  \label{eqcat-rel-fonc1}
F\text{-}\mathrm{Isoc} ^{\dag \dag} ( Y,\, (Y _\alpha, \PP , T_\alpha,  X) _{\alpha \in \Lambda}/K)
\rightarrow
F\text{-}\mathrm{Isoc} ^{\dag} ( Y'/K)
\times _{F\text{-}\mathrm{Isoc}  ( Y'/K)}
F\text{-}\mathrm{Isoc}  ( Y/K).
\end{equation}
\normalsize
Or, on bénéficie de l'équivalence de catégories (voir la preuve de \cite[6.5.3]{caro_devissge_surcoh}) :
\begin{equation}
  \label{eqcat-rel-fonc2}
F\text{-}\mathrm{Isoc} ^{\dag} ( Y'/K)
\times _{F\text{-}\mathrm{Isoc}  ( Y'/K)}
F\text{-}\mathrm{Isoc}  ( Y/K)
\cong
F\text{-}\mathrm{Isoc} ^\dag ( Y/K).
\end{equation}
On obtient alors
$F\text{-}\mathrm{Isoc} ^{\dag \dag} ( Y,\, (Y _\alpha, \PP , T_\alpha,  X) _{\alpha \in \Lambda}/K)
\rightarrow
F\text{-}\mathrm{Isoc} ^\dag ( Y/K)$,
que l'on notera $\rho _Y$,
en composant \ref{eqcat-rel-fonc1} et \ref{eqcat-rel-fonc2}.

$(b)$ Soient
$( \E _\alpha, \vartheta _{\alpha \beta})\in
F\text{-}\mathrm{Isoc} ^{\dag \dag} ( Y,\, (Y _\alpha, \PP , T_\alpha,  X) _{\alpha \in \Lambda}/K)$.
Avec l'étape $(a)$, on pose $E:= \rho _Y ( \E _\alpha, \vartheta _{\alpha \beta})$.
Construisons alors
$\E \in
F\text{-}\mathrm{Isoc} ^{\dag \dag} (\PP, \, T , \, X/K)$ tel que
$\mathcal{L} oc (\E) \riso ( \E _\alpha, \vartheta _{\alpha \beta})$.
Or,
$\sp _{Y _\alpha,+}( j^\dag _{T _\alpha} E)  | \U _\alpha \riso
\sp _{Y _\alpha \hookrightarrow \U _\alpha, +} ( \widehat{E} |\U _{\alpha K})
\riso$ \linebreak[3]
$\E _\alpha  | \U _\alpha  $ (voir \cite[5.2.4]{caro_devissge_surcoh}).
Par pleine fidélité du foncteur restriction $|\U _\alpha $ (voir \cite[6.5.2]{caro_devissge_surcoh}), il en découle
l'isomorphisme $\sp _{Y _\alpha,+}( j^\dag _{T _\alpha} E) \riso \E _\alpha$.
D'après \ref{lemm-eqcat-rel}, cela implique que
$\E _\alpha $
est $\D ^\dag _{\PP} (\hdag T) _{\Q}$-surcohérent.
En notant $\E _{\alpha \beta} := (\hdag T _\beta) (\E _\alpha)$, on dispose des deux flèches
$\E _\alpha \rightarrow \E _{\alpha \beta}$ et
$\E _\beta \rightarrow \E _{\beta \alpha} \underset{\vartheta _{\alpha \beta}}{\riso} \E _{\alpha \beta}$.
On pose alors
\begin{equation}
\notag
  \E := \ker ( \xymatrix { { \oplus _{\alpha \in \Lambda} \E _\alpha} \ar@<1ex>[r]\ar@<-1ex>[r]
  &
  {\oplus _{\alpha, \beta \in \Lambda}  \E _{\alpha \beta}} }).
\end{equation}
Comme pour tous $\alpha , \beta \in \Lambda$,
$\E _\alpha , \E _{\alpha \beta}$ sont
$\D ^\dag _{\PP} (\hdag T) _{\Q}$-surcohérents,
$\E$ est $\D ^\dag _{\PP} (\hdag T) _{\Q}$-surcohérent
(voir les propriétés de stabilité de la surcohérence de \cite{caro_surcoherent}).
Comme, pour tout $\alpha\in \Lambda$,
$ \E (\hdag T _\alpha) \riso \E _\alpha \in F\text{-}\mathrm{Isoc} ^{\dag \dag} (\PP , T _\alpha , X /K)$,
d'après \ref{isocsurcohgen-prel}, on obtient même $\E \in F\text{-}\mathrm{Isoc} ^{\dag \dag} (\PP , T , X /K)$.
Par construction, $\mathcal{L} oc (\E) \riso ( \E _\alpha, \vartheta _{\alpha \beta})$.
\hfill \hfill \qed \end{proof}

\begin{prop}
\label{ind-isosurcoh}
Avec les notations et hypothèses de \ref{isocsurcvgen-defi},
la catégorie
$F\text{-}\mathrm{Isoc} ^{\dag \dag}
  (Y,\, (Y _\alpha, \PP _\alpha, T_\alpha,  X_\alpha) _{\alpha \in \Lambda})$
  ne dépend pas du choix de la famille
$(Y _\alpha, \PP _\alpha, T_\alpha,  X_\alpha) _{\alpha \in \Lambda}$.
On la notera alors $F\text{-}\mathrm{Isoc} ^{\dag \dag}( Y/K )$.

Plus généralement, si $Y$ est un $k$-schéma séparé, lisse et $(Y _{l}) _{l=1,\dots ,N}$ sont ses composantes connexes,
on pose
$F\text{-}\mathrm{Isoc} ^{\dag \dag} (Y /K):= \prod _l F\text{-}\mathrm{Isoc} ^{\dag \dag } (Y _l/K)$.
 Ses objets sont les
{\it $F$-isocristaux surcohérents sur $Y$}.
\end{prop}

\begin{proof}
Il suffit de vérifier la validité des trois étapes analogues à celles de la preuve de \ref{ind-isocsurcvgen-vide}.
Pour l'étape i), cela correspond à \ref{eqcat-rel}. L'étape ii) s'établit de façon analogue tandis que pour
la dernière, il s'agit de remplacer l'utilisation de l'équivalence de catégories
\cite[2.3.5.i)]{Berig} par celle de
\cite[6.4.3.a)]{caro_devissge_surcoh}.
\hfill \hfill \qed \end{proof}

\begin{rema}
\label{rhoygen}
Grâce à \cite[6.3.3]{caro_devissge_surcoh}, la définition de \ref{ind-isosurcoh}
est compatible avec celle que l'on connaissait (i.e., \cite[6.4.3.b)]{caro_devissge_surcoh}).
\end{rema}

\subsection{\'Equivalence entre $F$-isocristaux surconvergents et surcohérents}

\begin{theo}
\label{eqcat-gen}
Soit $Y $ une $k$-variété lisse. On dispose d'une équivalence canonique de catégories
$$\sp _{Y,+}\ :\ F\text{-}\mathrm{Isoc} ^{\dag}( Y/K )\cong F\text{-}\mathrm{Isoc} ^{\dag \dag}( Y/K ).$$
\end{theo}
\begin{proof}
Il ne coûte rien de supposer $Y$ intègre.
  En choisissant un recouvrement de $Y$ satisfaisant les conditions de
  \ref{nota-Fisocdagetc}, cela résulte de \ref{congcasaffliss} (et de \ref{ind-isocsurcvgen-vide}, \ref{ind-isosurcoh}).
\hfill \hfill \qed \end{proof}

Il résulte de cette équivalence de catégories la description suivante des $F$-isocristaux surcohérents :

\begin{theo}
\label{surcoh=dagdag}
Avec les notations \ref{nota-isoc-surcoh}, la catégorie
$  F\text{-}\mathrm{Isoc} ^{\dag \dag}( Y/K)$ est équivalente à celle des
$F\text{-}\D ^\dag _{\PP} (\hdag T) _\Q$-modules surcohérents $\E$ à support dans $X$ tels que
$\E | \U$ soit dans l'image essentielle de $\sp _{Y \hookrightarrow \U, +}$.
\end{theo}

\begin{proof}
Soit $\E$ un $F\text{-}\D ^\dag _{\PP} (\hdag T) _\Q$-module surcohérent à support dans $X$ tels que
$\E | \U$ soit dans l'image essentielle de $\sp _{Y \hookrightarrow \U, +}$.
En notant $X_1 ,\dots , X _r$ les composantes irréductibles de $X$,
le morphisme canonique
$\oplus _r \R \underline{\Gamma} ^\dag  _{X _r} (\E) \rightarrow \E$
est un isomorphisme car il l'est en dehors de $T$ (voir \cite[4.3.12]{Be1}).
Par \cite[6.3.3]{caro_devissge_surcoh}, on peut alors supposer $X$ irréductible.
Avec ses notations et hypothèses, on dispose du diagramme \ref{diag631dejong}.
  Via la caractérisation de l'image essentielle de
  $\sp _{X' \hookrightarrow \PP',T',+}$ (voir \cite[4.1.9]{caro-construction}),
  il existe
  un (unique) $F$-isocristal $E'$  surconvergent sur $ Y '$
  tels que
$\R \underline{\Gamma} ^\dag _{X'} f_T ^! (\E ) \riso \sp _{X' \hookrightarrow \PP',T',+}(E') $.
Grâce à l'équivalence de catégories \ref{eqcat-rel-fonc2},
il existe donc un unique (à isomorphisme près) $F$-isocristal surconvergent $E$ sur $Y$ induisant $E'$ et
le $F$-isocristal convergent sur $Y$ correspondant à $\E |\U$.

Avec \ref{eqcat-gen}, soit
$\widetilde{\E} := \sp _{Y,+} (E) \in F\text{-}\mathrm{Isoc} ^{\dag \dag}( \PP, T, X/K)$.
En utilisant le théorème de pleine fidélité du foncteur restriction en dehors de $T'$
(voir \cite[6.5.2]{caro_devissge_surcoh}),
on vérifie $ \R \underline{\Gamma} ^\dag _{X'} f_T ^! (\widetilde{\E} ) \riso \sp _{X' \hookrightarrow \PP',T',+}(E')$.
D'où $ \R \underline{\Gamma} ^\dag _{X'} f_T ^! (\widetilde{\E} ) \riso$ \linebreak[3]
$ \R \underline{\Gamma} ^\dag _{X'} f_T ^! (\E)$.
Or, lors de la preuve de \cite[6.3.1]{caro_devissge_surcoh},
on a construit un morphisme
$ \widetilde{\E} \rightarrow f _{T,+}  \R \underline{\Gamma} ^\dag _{X'} f_T ^! (\widetilde{\E} )$.
En le composant avec le morphisme d'adjonction
$f _{T,+} \R \underline{\Gamma} ^\dag _{X'} f_T ^! (\E ) \rightarrow \E$, on obtient
$\widetilde{\E} \rightarrow \E$ qui est par construction un isomorphisme en dehors de $T$.
On conclut la preuve avec \cite[4.3.12]{Be1}.
\hfill \hfill \qed \end{proof}

\section{Autour de la dévissabilité en $F$-isocristaux surconvergents}

Soit $\PP$ un $\V$-schéma formel propre et lisse.
Si $Y$ est un sous-schéma de $P$, on notera $\overline{Y}$ l'adhérence schématique de $Y$ dans $P$.

\subsection{Dévissage des $F$-complexes surcohérents en $F$-isocristaux surconvergents}

Déduisons maintenant du théorème \ref{surcoh=dagdag} le lemme et le théorème suivants.

\begin{lemm}
\label{lemm-theodev2}
  Soient $T$ un diviseur de $P$, $Y $ un sous-schéma fermé de $P \setminus T $,
  $Y _1$ une composante irréductible de $Y$ et
  $\E $ un $F$-$\D ^\dag _{\PP} (\hdag T) _\Q$-module surcohérent à support dans $\overline{Y}$.

  Il existe un diviseur $T '\supset T$ de $P$ tel que,
  en notant $Y' := \overline{Y} \setminus T'$,
  $Y'$ soit intègre, lisse, inclus et dense dans $Y _1$,
  $\E (\hdag T') \in F\text{-}\mathrm{Isoc} ^{\dag \dag}( \PP, T', \overline{Y}/K)=
  F\text{-}\mathrm{Isoc} ^{\dag \dag}( Y'/K)$.
\end{lemm}

\begin{proof}
Comme $P$ est somme de ses composantes irréductibles, il ne coûte rien de supposer
$P$ intègre.
Soit $\U'$ un ouvert affine de $\PP$ inclus dans $\PP \setminus T$
tel que $U '\cap Y$ soit intègre, lisse, inclus et dense dans $Y _1$.
Ainsi, $\E | \U'$ est à support dans $U'\cap Y $.
Via \cite[2.2.17]{caro_courbe-nouveau}, quitte à diminuer $\U'$, on peut même supposer
que $\E | \U'$
est dans l'image essentielle de $\sp _{U'\cap Y \hookrightarrow \U',+}$ (voir \cite[4.1.9]{caro-construction}).
Avec \ref{affdiv},
$T':= P \setminus U'$ est un diviseur. Le théorème \ref{surcoh=dagdag} nous permet de conclure.
\hfill \hfill \qed \end{proof}

\begin{theo}
\label{caradev}
Soient $T$ un diviseur de $P$, $Y $ un sous-schéma fermé de $P \setminus T $,
  $\E \in F\text{-}D ^\mathrm{b} _\mathrm{surcoh} (\smash{\D} ^\dag _{\PP} (\hdag T) _\Q)$
  à support dans $\overline{Y}$.

Il existe alors des diviseurs $T _1, \dots,T _{r}$ contenant $T$ avec $T _r =T$ tels que,
en notant $T _0: =\overline{Y}$,
pour $i=0,\dots , r-1$,
  le schéma $Y _i:= T _0 \cap \dots \cap T _i \setminus T _{i+1}$
  soit lisse et, pour tout entier $j$,
\begin{gather}
  \notag
  \mathcal{H} ^j (\R \underline{\Gamma} ^\dag _{T _0 \cap \dots \cap T _i}  (\hdag T _{i+1}) (\E) )
  \in F\text{-}\mathrm{Isoc} ^{\dag \dag}( \PP, T _{i+1}, T _0 \cap \dots \cap T _i/K)
  \\ \notag
  =F\text{-}\mathrm{Isoc} ^{\dag \dag}( Y_i/K) \cong
  F\text{-}\mathrm{Isoc} ^{\dag}( Y_i/K).
\end{gather}

\end{theo}

\begin{proof}
On utilise successivement \ref{lemm-theodev2} jusqu'à ce que l'ouvert dense soit l'espace tout entier
(i.e., avec les notations de \ref{lemm-theodev2}, $T' =T$).
\hfill \hfill \qed \end{proof}

\begin{rema}
\label{remacaradev}
 De la même manière, on établit que le lemme \ref{lemm-theodev2} est encore valable en remplaçant
  {\og $Y'$ soit intègre, lisse, inclus et dense dans $Y _1$ \fg}
  par
  {\og $Y'$ soit intègre, affine, lisse, inclus et dense dans $Y _1$ et muni d'un modèle idéal\fg}
  (voir \cite[7.4.1 et 7.4.2]{caro_devissge_surcoh}).
  On vérifie alors que le théorème \ref{caradev} reste valable en remplaçant
   {\og le schéma $Y _i:=Y \cap T _0 \cap \dots \cap T _i \setminus T _{i+1}$
  soit lisse\fg} par
     {\og le schéma $Y _i:=Y \cap T _0 \cap \dots \cap T _i \setminus T _{i+1}$
  soit affine, lisse et muni d'un modèle idéal\fg}.
\end{rema}

\subsection{Sur la définition de la dévissabilité en $F$-isocristaux surconvergents}

\begin{vide}
  Soient $X$, $X'$, $T$, $T'$ des sous-schémas fermés de $P$ tels que
  $X \setminus T = X' \setminus T'$.
  Pour tout $\E \in (F\text{-})\smash{\underset{^{\longrightarrow}}{LD}} ^{\mathrm{b}} _{\Q ,\mathrm{qc}}
(\smash{\widehat{\D}} _{\PP} ^{(\bullet)})$,
on dispose de l'isomorphisme canonique :
\begin{equation}
  \label{xtx't'}
\R\underline{\Gamma} ^\dag _{X} (\hdag T ) (\E)
\riso
\R\underline{\Gamma} ^\dag _{X'} (\hdag T ') (\E).
\end{equation}
En effet, $\R\underline{\Gamma} ^\dag _{X} (\hdag T ) (\E) \riso
\R\underline{\Gamma} ^\dag _{X} (\hdag T ) (\O _{\PP,\Q})
\smash{\overset{\L}{\otimes}}   ^{\dag} _{\O  _{\PP,\Q}} \E$ (voir \cite[2.2.6.1]{caro_surcoherent}),
de même avec des primes. On se ramène ainsi au cas où $\E = \O _{\PP,\Q}$.
Or, comme $\O _{\PP, \Q}$ est
$\D ^\dag _{\PP,\Q}$-surcohérent, on déduit de \cite[2.2.8, 2.2.14, 3.2.4]{caro_surcoherent} l'isomorphisme
$\R\underline{\Gamma} ^\dag _{X} (\hdag T ) (\O _{\PP, \Q})
\riso
\R\underline{\Gamma} ^\dag _{X \cap X'} (\hdag T \cup T ') (\O _{\PP, \Q})$.
D'où le résultat par symétrie.

\noindent $\bullet$ En posant $ Y := X \setminus T$, on note alors sans ambiguïté
$\R\underline{\Gamma} ^\dag _{Y} (\E) $ l'un
des deux ($F$-)complexes de \ref{xtx't'}.

\noindent $\bullet$ Si $Y $ et $Y'$ sont deux sous-schémas de $P$, pour tout
$\E \in (F\text{-})\smash{\underset{^{\longrightarrow}}{LD}} ^{\mathrm{b}} _{\Q ,\mathrm{qc}}
(\smash{\widehat{\D}} _{\PP} ^{(\bullet)})$,
il résulte de \cite[2.2.8 et 2.2.14]{caro_surcoherent} l'isomorphisme canonique :
\begin{equation}
  \label{gammayY'}
  \R\underline{\Gamma} ^\dag _{Y} \circ \R\underline{\Gamma} ^\dag _{Y'} (\E)
  \riso
  \R\underline{\Gamma} ^\dag _{Y \cap Y'} (\E).
\end{equation}

\noindent $\bullet$ Soit $\E \in (F\text{-})\smash{\underset{^{\longrightarrow}}{LD}} ^{\mathrm{b}} _{\Q ,\mathrm{qc}}
(\smash{\widehat{\D}} _{\PP} ^{(\bullet)})$.
Si $Y '$ est un ouvert (resp. un fermé) de $Y$, on dispose du morphisme canonique
$\R\underline{\Gamma} ^\dag _{Y} (\E) \rightarrow \R\underline{\Gamma} ^\dag _{Y'} (\E)$
(resp. $\R\underline{\Gamma} ^\dag _{Y'} (\E) \rightarrow \R\underline{\Gamma} ^\dag _{Y} (\E)$).
Si $Y'$ est un fermé de $Y$, on bénéficie aussi du triangle distingué de localisation :
$\R\underline{\Gamma} ^\dag _{Y'} (\E) \rightarrow \R\underline{\Gamma} ^\dag _{Y} (\E)
\rightarrow \R\underline{\Gamma} ^\dag _{Y \setminus Y'} (\E) \rightarrow +1.$

\end{vide}

\begin{vide}
\label{def-espcohoisoc}
  Soient $X$ et $X'$ deux sous-schémas fermés de $P$, $T$ et $T'$ deux diviseurs de $P$ tels que
  $X \setminus T = X' \setminus T'$. On suppose $Y := X \setminus T$ lisse. Pour tous
  $\E \in F\text{-}\smash{\underset{^{\longrightarrow}}{LD}} ^{\mathrm{b}} _{\Q ,\mathrm{qc}}
(\smash{\widehat{\D}} _{\PP} ^{(\bullet)})$, $j \in \Z$,
$\mathcal{H} ^j (\underset{\longrightarrow}{\lim} \R\underline{\Gamma} ^\dag _{X} (\hdag T ) (\E) )
\in F \text{-} \mathrm{Isoc} ^{\dag \dag} ( \PP, T ,X /K)$
si et seulement si
$\mathcal{H} ^j (\underset{\longrightarrow}{\lim} \R\underline{\Gamma} ^\dag _{X'} (\hdag T ') (\E) )
\in F \text{-} \mathrm{Isoc} ^{\dag \dag} ( \PP, T ',X '/K)$ (cela découle de \cite[6.3.4]{caro_devissge_surcoh}).

On dira alors sans ambiguïté que {\og les espaces de cohomologie de
  $\R \underline{\Gamma} ^\dag _{Y}  (\E)$
  sont associés à des $F$-isocristaux surconvergents sur $Y $ \fg}
  si, pour tout entier $j$, $\mathcal{H} ^j (\underset{\longrightarrow}{\lim} \R \underline{\Gamma} ^\dag _{Y}  (\E))
  \in F \text{-} \mathrm{Isoc} ^{\dag \dag} ( Y/K) \cong
  F \text{-} \mathrm{Isoc} ^{\dag} ( Y/K) $ (voir \ref{eqcat-gen}).
\end{vide}
Le lemme suivant est aisé.
\begin{lemm}
  \label{def-espcohoisoc-tri}
  Soit $Y$ un sous-schéma lisse de $P$ tel qu'il existe un diviseur $T$ de $P$
  vérifiant $Y=\overline{Y} \setminus T $.
  La sous-catégorie pleine de $F\text{-}\smash{\underset{^{\longrightarrow}}{LD}} ^{\mathrm{b}} _{\Q ,\mathrm{qc}}
(\smash{\widehat{\D}} _{\PP} ^{(\bullet)})$
des $F$-complexes $\E$
tels que les espaces de cohomologie de
  $\R \underline{\Gamma} ^\dag _{Y}  (\E)$
soient associés à des $F$-isocristaux surconvergents sur $Y $ est triangulée.
\end{lemm}

\begin{lemm}
\label{stab-espcohoisoc}
  Soient $Y$, $Y'$ deux sous-schémas de $P$ tels qu'il existe des diviseurs $T$, $T'$ de $P$
  satisfaisant $Y = \smash{\overline{Y}} \setminus T$ et $Y '= \smash{\overline{Y}} '\setminus T'$.
  On suppose $Y$ et $Y \cap Y'$ lisses.
  Soit $\E \in F \text{-} \smash[b]{\underset{^{\longrightarrow }}{LD }}  ^\mathrm{b} _{\Q, \mathrm{qc}}
(\overset{^\mathrm{g}}{} \smash{\widehat{\D}} _{\PP} ^{(\bullet)})$ tel que
les espaces de cohomologie de
  $\R \underline{\Gamma} ^\dag _{Y}  (\E)$
  soient associés à des $F$-isocristaux surconvergents sur $Y $. Alors les espaces de cohomologie de
$\R \underline{\Gamma} ^\dag _{Y\cap Y'}  (\E)$
  sont associés à des $F$-isocristaux surconvergents sur $Y \cap Y'$.
\end{lemm}
\begin{proof}
  Cela dérive de \cite[6.3.2 et 6.3.5]{caro_devissge_surcoh} et de \ref{gammayY'}.
\hfill \hfill \qed \end{proof}

\'Etendons à présent
la notion de dévissabilité en $F$-isocristaux surconvergents définie en \cite[8.1.1]{caro_devissge_surcoh}
au cas des $F$-complexes quasi-cohérents. Une première définition est donnée ci-dessous, une seconde plus
générale est obtenue en \ref{coro-caradev2}.

\begin{defi}
\label{coro-caradev}
Soient $\E \in F \text{-} \smash[b]{\underset{^{\longrightarrow }}{LD }}  ^\mathrm{b} _{\Q, \mathrm{qc}}
(\overset{^\mathrm{g}}{} \smash{\widehat{\D}} _{\PP} ^{(\bullet)})$,
$Y$ un sous-schéma de $P$.
On suppose qu'il existe un diviseur $T$ de $P$ tel que $Y =\smash{\overline{Y}} \setminus T$.

Le $F$-complexe $\E$ {\og se dévisse au dessus de $Y$ en $F$-isocristaux surconvergents\fg}
s'il existe des diviseurs $T _1, \dots,T _{r}$ contenant $T$ avec $T _r =T$ tels que,
pour $i: =0,\dots, r-1$ et avec $T _0:=\smash{\overline{Y}}$,
$Y _i := T _0 \cap T _1 \cap \dots \cap T _i \setminus T _{i+1}$
  soit lisse et
  les espaces de cohomologie de
  $\R \underline{\Gamma} ^\dag _{Y _i} (\E)$
  soient associés à des $F$-isocristaux surconvergents sur $Y _i$ (voir \ref{def-espcohoisoc}).

On dira aussi que le $F$-complexe {\og $\E$ se dévisse au-dessus de
la stratification
$Y =\sqcup _{i= 0,\dots , r-1} Y _i$ en $F$-isocristaux surconvergents\fg}
ou que {\og $(T _1, \dots, T_r)$ fournit un dévissage sur $Y$ de $\E$ en $F$-isocristaux surconvergents\fg}.
\end{defi}

\begin{rema}
\label{rema-caradev}
Avec les notations et hypothèses de \ref{coro-caradev},
le fait qu'un $F$-complexe $\E$ de $F \text{-} \smash[b]{\underset{^{\longrightarrow }}{LD }}  ^\mathrm{b} _{\Q, \mathrm{qc}}
(\overset{^\mathrm{g}}{} \smash{\widehat{\D}} _{\PP} ^{(\bullet)})$
 {\og se dévisse au-dessus de
la stratification
$Y =\sqcup _{i= 0,\dots , r-1} Y _i$
en $F$-isocristaux surconvergents \fg}
ne dépend ni du choix du diviseur $T$ tel que $Y =\smash{\overline{Y}} \setminus T$, ni de celui
du $r$-uplet $(T _1, \dots, T _r)$ où
les $T _i$ sont des diviseurs contenant $T _r$ tels que
 $Y _i =  T _0 \cap \dots \cap T _i \setminus T _{i+1}$
où $T _0:=\smash{\overline{Y}}$ (on remarque que cela implique $Y= \smash{\overline{Y}} \setminus T _r$).
Les définitions de \ref{coro-caradev} n'engendrent donc aucune ambiguïté.
On remarque enfin que cela ne change rien si on remplace
$\smash{\overline{Y}}$ par un sous-schéma fermé $X$ de $P$ tel que $Y=X \setminus T $.
\end{rema}

\begin{prop}
\label{devdev=dev}
Soient $\E \in F \text{-} \smash[b]{\underset{^{\longrightarrow }}{LD }}  ^\mathrm{b} _{\Q, \mathrm{qc}}
(\overset{^\mathrm{g}}{} \smash{\widehat{\D}} _{\PP} ^{(\bullet)})$,
$Y$ un sous-schéma de $P$ tel qu'il existe un diviseur $T$ de $P$
satisfaisant $Y =\smash{\overline{Y}} \setminus T$.
Soient $T _1,\dots, T _{r}$ des diviseurs de $P$ contenant $T$ avec $T _r =T$ et,
pour $i: =0,\dots, r-1$,
$Y _i := T _0 \cap T _1 \cap \dots \cap T _i \setminus T _{i+1}$ avec $T _0:=\smash{\overline{Y}}$.

Si, pour tout $i: =0,\dots, r-1$, $\E$
est dévissable sur $Y _i$ en $F$-isocristaux surconvergents alors
$\E$ l'est sur $Y$.
\end{prop}

\begin{proof}
Plus précisément, nous prouvons le lemme \ref{lemm-devdev=dev} ci-après :
\begin{lemm}
  \label{lemm-devdev=dev}

On garde les notations et hypothèses de \ref{devdev=dev}.
Pour tout $i =0,\dots , r-1$,
il existe donc
des diviseurs $T _{(i, 1)}, \dots,T _{(i, r _i)}$ contenant $T _{i+1}$ avec $T _{(i, r _i)}=T _{i+1}$ tels que,
en notant $T _{(i, 0)} :=\overline{Y} _i $,
pour $h=0,\dots , r_i-1$,
  le $k$-schéma $Y _{(i,h)}:= T_{(i,0)} \cap \dots \cap T _{(i,h)} \setminus T _{(i,h+1)}$
  soit lisse et, pour tout entier $j$,
  $$ \mathcal{H} ^j (\underset{\longrightarrow}{\lim}
  \R \underline{\Gamma} ^\dag _{Y _{(i,h)}}  \E  )
  \in F\text{-}\mathrm{Isoc} ^{\dag \dag}(Y _{(i,h)}/K)
  \cong F\text{-}\mathrm{Isoc} ^{\dag}( Y_{(i,h)}/K).$$

Alors $(T _{(0, 1)}, \dots,T _{(0, r _0)},
T _{(1, 1)}, \dots,T _{(1, r _1)},\dots ,
T _{(r-1, 1)}, \dots,T _{(r-1, r _{r-1} )})$
fournit un dévissage en $F$-isocristaux surconvergents de $\E$
au-dessus de la stratification
\small
\begin{equation}
\label{strat-comp}
  Y = Y _{(0,0)} \sqcup \dots \sqcup Y _{(0,r _0 -1)} \bigsqcup
Y _{(1,0)} \sqcup \dots \sqcup Y _{(1,r _1 -1)} \bigsqcup
\dots \bigsqcup Y _{(r -1,0)} \sqcup \dots \sqcup Y _{(r- 1,r _{r-1} -1)}.
\end{equation}
\normalsize
\end{lemm}

\begin{proof}
  On procède par récurrence sur $r$. Lorsque $r =1$, c'est tautologique.
Supposons donc la proposition valable pour $r-1\geq 1$ et prouvons-le pour $r$.
On dispose de la décomposition :
\scriptsize
\begin{equation}
\label{devdev=dev-=}
  Y = (Y \setminus T _{(0,1)} )\bigsqcup (Y \cap T _{(0,1)} \setminus T _{(0,2)})
\bigsqcup  \dots \bigsqcup
(Y \cap T _{(0,1)} \cap \dots \cap T _{(0,r _0 -1)} \setminus T _{(0, r_0)})
\bigsqcup
(Y \cap T _{(0,1)} \cap \dots \cap T _{(0,r _0)}).
\end{equation}
\normalsize
Pour tout $h=0,\dots ,r_0 -1$, comme $T _{(0,h)}\supset T _1 \supset T$
(resp. et $T _{(0, r _0)}=T _{1}$),
on obtient
$Y \cap T _{(0,1)} \cap \dots \cap T _{(0,h )} \setminus T _{(0, h+1)} = Y _{(0,h)} $
(resp. $Y \cap T _{(0,1)} \cap \dots \cap T _{(0,r _0)}
=
Y \cap T _1$).
La décomposition \ref{devdev=dev-=} est ainsi exactement :
\begin{equation}
\label{devdev=dev-=2}
  Y = Y _{(0,0)}\sqcup Y _{(0,1)}
\sqcup  \dots \sqcup
Y _{(0,r _0 -1)}
\sqcup
(Y \cap T _1).
\end{equation}
La décomposition \ref{devdev=dev-=2} signifie donc
que $(T _{(0, 1)}, \dots,T _{(0, r _0)})$ fournit le début d'un dévissage sur $Y$
(le dévissage ne s'arrête pas car à priori $T _{(0, r _0)} \not =T$, i.e. $Y \cap T _1 \not = \emptyset$)
de $\E$ en $F$-isocristaux surconvergents.
Enfin, par hypothèse de récurrence appliquée à $\E$ sur $Y \cap T _1$ pour $(T _2, \dots, T _r)$,
la famille de diviseurs $(T _{(1, 1)}, \dots,T _{(1, r _1)},\dots ,
T _{(r-1, 1)}, \dots,T _{(r-1, r _{r-1} )})$
induit un dévissage en $F$-isocristaux surconvergents de
$\E$
au-dessus de la stratification $Y \cap T _1= Y _{(1,0)} \sqcup \dots \sqcup Y _{(1,r _1 -1)} \bigsqcup
\dots \bigsqcup Y _{(r -1,0)} \sqcup \dots \sqcup Y _{(r- 1,r _{r-1} -1)}$.
\hfill \hfill \qed \hfill \hfill \qed \end{proof}

\end{proof}

\begin{lemm}
\label{indYdev-lemm1}
Soient $\E \in F \text{-} \smash[b]{\underset{^{\longrightarrow }}{LD }}  ^\mathrm{b} _{\Q, \mathrm{qc}}
(\overset{^\mathrm{g}}{} \smash{\widehat{\D}} _{\PP} ^{(\bullet)})$,
$Y ' \subset Y$ deux sous-schémas de $P$.
On suppose qu'il existe deux diviseurs $T$, $T'$ de $P$
tels que
$Y ={\overline{Y}} \setminus T$ et $Y' ={\overline{Y}} '\setminus T'$.
Si $\E$ est dévissable sur $Y$ en $F$-isocristaux surconvergents
alors $\E$ l'est sur $Y'$.
\end{lemm}

\begin{proof}
Quitte à remplacer $T'$ par $T \cup T '$, on peut supposer $T ' \supset T$.
    Soient $T _1, \dots,T _{r}$ des diviseurs de $P$ contenant $T$ avec $T _r =T$ tels que
    $(T _1,\dots , T _r)$ fournisse un dévissage sur $Y$ de $\E$ en $F$-isocristaux surconvergents.
Notons $T _0 :=\smash{\overline{Y}}$, $T ' _0: =\smash{\overline{Y'}}$ et, pour $i=0,\dots , r-1$,
$T ' _{i+1}:= T _{i+1} \cup T'$,
$Y  _i:=T _0 \cap \dots \cap T _i \setminus T _{i+1}$, $Y ' _i:=T '_0 \cap \dots \cap T '_i \setminus T '_{i+1}$.
Pour $i=0,\dots , r-1$, comme $Y ' _i =  Y _i \cap \smash{\overline{Y}}' \setminus T ' _{i+1}$ et
comme $\R \underline{\Gamma} ^\dag _{Y  _i} \E$ est $\D ^\dag _{\PP} (\hdag T _{i+1}) _{\Q}$-surcohérent,
par stabilité de la surcohérence par foncteur cohomologique local à support strict dans un sous-schéma fermé
(voir \cite[3.1.5]{caro_surcoherent}),
il en résulte que
$\R \underline{\Gamma} ^\dag _{Y ' _i} \E$ est $\D ^\dag _{\PP} (\hdag T '_{i+1}) _{\Q}$-surcohérent.
Par \ref{caradev},
il en découle que $\E$ se dévisse sur $Y '_i $ en $F$-isocristaux surconvergents.
La proposition \ref{devdev=dev} nous permet de conclure.
\hfill \hfill \qed \end{proof}

\begin{lemm}
\label{EE'dev}
Soient $\E, \E ' \in F \text{-} \smash[b]{\underset{^{\longrightarrow }}{LD }}  ^\mathrm{b} _{\Q, \mathrm{qc}}
(\overset{^\mathrm{g}}{} \smash{\widehat{\D}} _{\PP} ^{(\bullet)})$,
$Y$ un sous-schéma de $P$ tel qu'il existe un diviseur $T$ de $P$ vérifiant
$Y =\smash{\overline{Y}} \setminus T$.
Si $\E$ et $\E'$ se dévissent sur $Y$ en $F$-isocristaux surconvergents alors,
il existe des diviseurs $T _1, \dots,T _{r}$ contenant $T$ avec $T _r =T$ tels que
$(T _1, \dots,T _{r})$ induise un dévissage sur $Y$ de $\E$ et $\E'$ en $F$-isocristaux surconvergents.

La sous-catégorie pleine de $F \text{-} \smash[b]{\underset{^{\longrightarrow }}{LD }}  ^\mathrm{b} _{\Q, \mathrm{qc}}
(\overset{^\mathrm{g}}{} \smash{\widehat{\D}} _{\PP} ^{(\bullet)})$
des $F$-complexes qui se dévissent sur $Y$ en $F$-isocristaux surconvergents est alors triangulée.
\end{lemm}

\begin{proof}
Il existe des diviseurs $T _1, \dots,T _{r}$ contenant $T$ avec $T _r =T$ tels que
$(T _1, \dots,T _{r})$ donne un dévissage sur $Y$ de $\E$ en $F$-isocristaux surconvergents.
Notons $T _0:=\smash{\overline{Y}}$ et, pour $i =0,\dots, r-1$,
$Y _i := T _0 \cap \dots \cap T _i \setminus T _{i+1}$.
Il découle de \ref{indYdev-lemm1} que $\E'$ se dévisse sur $Y _i$ en $F$-isocristaux surconvergents.
D'après \ref{lemm-devdev=dev} et avec notations,
$\E'$ se dévisse en $F$-isocristaux surconvergents
au-dessus d'une stratification de la forme \ref{strat-comp}.
En utilisant \ref{stab-espcohoisoc}, on vérifie que $\E$ l'est de même sur cette stratification \ref{strat-comp}.

La dernière assertion résulte de \ref{def-espcohoisoc-tri}.
\hfill \hfill \qed \end{proof}

\begin{lemm}
\label{EdagtTT'}
  Soient $T, T '$ des diviseurs de $P$, $X$ un sous-schéma fermé de $P$,
    $\E \in F \text{-} \smash[b]{\underset{^{\longrightarrow }}{LD }}  ^\mathrm{b} _{\Q, \mathrm{qc}}
(\overset{^\mathrm{g}}{} \smash{\widehat{\D}} _{\PP} ^{(\bullet)})$.

\begin{enumerate}
  \item \label{EdagtTT'i} Si le $F$-complexe $\E$ se dévisse sur $X \setminus T'$ en $F$-isocristaux surconvergents alors
$\E (\hdag T)$ l'est sur $X \setminus (T \cup T')$.
  \item \label{EdagtTT'ii} Le $F$-complexe $\E (\hdag T)$ se dévisse sur $X \setminus T'$ en $F$-isocristaux surconvergents
  si et seulement si $\E (\hdag T)$ l'est sur $X \setminus (T \cup T')$.
\end{enumerate}
\end{lemm}

\begin{proof}
  Si $(T ' _1, \dots , T ' _r)$ fournit un dévissage de $\E$ sur $X \setminus T'$ en $F$-isocristaux surconvergents,
  alors
 $(T ' _1 \cup T, \dots , T ' _r \cup T)$ en fournit un sur $X \setminus (T \cup T')$ à $\E (\hdag T)$.
 D'où \ref{EdagtTT'i}).
 La nécessité de \ref{EdagtTT'ii}) découle de \ref{indYdev-lemm1} (ou de \ref{EdagtTT'i})).
 Réciproquement,
si $\E (\hdag T)$ se dévisse sur $X \setminus (T '\cup T) $ en $F$-isocristaux surconvergents,
comme
 $\R \underline{\Gamma} ^\dag _{(X \setminus T')\cap T}  (\E (\hdag T))=0$
 et $X \setminus T'= (X \setminus (T '\cup T) )\sqcup ((X \setminus T')\cap T)$,
 la proposition \ref{devdev=dev} (on prend $T_1=T\cup T'$ et $T _2 =T'$) nous permet de conclure.
\hfill \hfill \qed \end{proof}

\begin{lemm}
\label{lemm-devloc}
  Soient $T', T _1,\dots , T _r$ des diviseurs de $P$, $X$ un sous-schéma fermé de $P$,
  $\E \in F \text{-} \smash[b]{\underset{^{\longrightarrow }}{LD }}  ^\mathrm{b} _{\Q, \mathrm{qc}}
(\overset{^\mathrm{g}}{} \smash{\widehat{\D}} _{\PP} ^{(\bullet)})$. On pose
  $T := \cap _{i=1} ^r T _i$, $Y ' := X \setminus T'$.

  Le $F$-complexe $\E (\hdag T)$ est dévissable sur $Y'$ en $F$-isocristaux surconvergents si et seulement si, pour tout
  $i=1, \dots , r$, $\E (\hdag T _i)$ l'est sur $Y'$.
\end{lemm}

\begin{proof}
D'après respectivement \ref{EdagtTT'}.\ref{EdagtTT'i}) et \ref{EdagtTT'ii}), si
$\E (\hdag T)$ est dévissable sur $Y'$ en $F$-isocristaux surconvergents
alors $\E (\hdag T _i)$ l'est sur $Y' \setminus T _i$ et donc sur $Y'$.
\'Etablissons
à présent la suffisance.
Pour cela, prouvons par récurrence sur $s $ tel que $ 1 \leq s \leq r$ que, pour tous diviseurs $T ' _1 ,\dots , T ' _s$
chacun égaux à des réunions de
$T_1, \dots, T _r$, $\E (\hdag T ' _1 \cap \dots \cap T ' _s)$ est dévissable sur $Y'$ en $F$-isocristaux surconvergents.
Le cas où $s =1$ découle de \ref{EdagtTT'}.\ref{EdagtTT'i}) et \ref{EdagtTT'ii}).
Soient $T ' _1 ,\dots , T ' _{s+1}$ des diviseurs chacun égaux à des réunions de
$T_1, \dots, T _r$ et, pour $i =1,\dots ,s$, $T '' _i:= T ' _i \cup T ' _{s+1}$.
On dispose du triangle distingué de Mayer-Vietoris (voir \cite[2.2.16]{caro_surcoherent}) :
\begin{equation}
  \notag
  \E (\hdag T ' _1 \cap \dots \cap T ' _{s+1}) \rightarrow
\E (\hdag T ' _1 \cap \dots \cap T ' _{s})
\oplus
\E (T ' _{s+1})
\rightarrow
\E (\hdag T '' _1 \cap \dots \cap T '' _{s})
\rightarrow +1 .
\end{equation}
Par hypothèse de récurrence, les termes du milieu et de droite sont dévissables sur $Y'$ en $F$-isocristaux surconvergents.
Avec \ref{EE'dev}, il en résulte qu'il en est de même du terme de gauche.
\hfill \hfill \qed \end{proof}

\begin{lemm}
  \label{lemgammaYdevup}
Soient $Y$ et $Y'$ des sous-schémas de $P$ tels qu'il existe des diviseurs $T$, $T'$ vérifiant
$Y = \smash{\overline{Y}}\setminus T$,
$Y' = \smash{\overline{Y}}'\setminus T'$.
 Soit   $\E \in F \text{-} \smash[b]{\underset{^{\longrightarrow }}{LD }}  ^\mathrm{b} _{\Q, \mathrm{qc}}
(\overset{^\mathrm{g}}{} \smash{\widehat{\D}} _{\PP} ^{(\bullet)})$ tel
que $\R \underline{\Gamma} ^\dag _{Y'} (\E) \riso
\R \underline{\Gamma} ^\dag _Y (\E)$. Alors $\E$ est dévissable sur $Y$ en $F$-isocristaux surconvergents
si et seulement s'il l'est sur $Y'$.
\end{lemm}

\begin{proof}

Quitte à remplacer $Y$ par $Y \cap Y'$ et $T$ par $T \cup T'$,
on se ramène via \ref{gammayY'} au cas où
$Y \subset Y'$ et $T \supset T'$.
La réciproque étant \ref{indYdev-lemm1}, supposons
$\E$ dévissable sur $Y$ en $F$-isocristaux surconvergents et
prouvons alors que $\E$ l'est sur $Y'$.
Or, en notant $\FF := \R \underline{\Gamma} ^\dag _{\smash{\overline{Y}}'} (\E)$,
$\E$ est dévissable sur $Y'$ en $F$-isocristaux surconvergents
si et seulement si $(\hdag T ') (\FF)$ l'est.
De plus, l'isomorphisme $\R \underline{\Gamma} ^\dag _{Y'} (\E) \riso
\R \underline{\Gamma} ^\dag _Y (\E)$ implique le suivant $(\hdag T ') (\FF) \riso (\hdag T) \FF$
(on utilise \cite[2.2.14]{caro_surcoherent}).
Par \ref{EdagtTT'}.\ref{EdagtTT'ii}), il suffit alors de prouver que $(\hdag T ') (\FF) $ est dévissable sur $Y '\setminus T$ en
 $F$-isocristaux surconvergents, ce qui revient à traiter le cas $T =T'$.
 Supposons donc $T =T'$.
  Soient $(T_1,\dots, T _r)$ des diviseurs fournissant sur $Y$ un dévissage de $\E$
  en $F$-isocristaux surconvergents.
Notons $T  _0: =\smash{\overline{Y}}'$ et, pour $i=0,\dots , r-1$,
$Y ' _i:=T _0 \cap \dots \cap T _i \setminus T _{i+1}$.
Par définition, les espaces de cohomologie de
 $\R \underline{\Gamma} ^\dag _{Y' _i \cap Y} (\E)$ sont associés à des $F$-isocristaux surconvergents sur
 $Y' _i \cap Y$.
 Or, avec \ref{gammayY'} et par hypothèse, on obtient
$\R \underline{\Gamma} ^\dag _{Y' _i } (\E)
=\R \underline{\Gamma} ^\dag _{Y' _i \cap Y'} (\E)
\riso
\R \underline{\Gamma} ^\dag _{Y' _i \cap Y} (\E)$.
Ainsi, $\R \underline{\Gamma} ^\dag _{Y' _i } (\E)$
est $\D ^\dag _{\PP} (\hdag T _{i+1}) _{\Q}$-surcohérent.
Par \ref{caradev},
cela implique que $\E$ se dévisse sur $Y '_i $ en $F$-isocristaux surconvergents.
Le $F$-complexe $\E$ l'est donc sur $Y'$ (voir \ref{devdev=dev}).
\hfill \hfill \qed \end{proof}

\begin{defi}
  \label{coro-caradev2}
Soient $\E \in F \text{-} \smash[b]{\underset{^{\longrightarrow }}{LD }}  ^\mathrm{b} _{\Q, \mathrm{qc}}
(\overset{^\mathrm{g}}{} \smash{\widehat{\D}} _{\PP} ^{(\bullet)})$,
$Y$ un sous-schéma de $P$.

Le $F$-complexe $\E$ est {\og dévissable sur $Y$ en $F$-isocristaux surconvergents\fg} s'il existe
un recouvrement fini ouvert $(Y _l) _l$ de $Y$,
des diviseurs $T _l$ de $P$ tels que, pour tout $l$,
$Y _l = \smash{\overline{Y}} _l \setminus T _l$ et $\E$ soit dévissable sur $Y _l$ en $F$-isocristaux surconvergents (voir \ref{coro-caradev}).
\end{defi}

\begin{prop}
  \label{devloc}
Soient $\E \in F \text{-} \smash[b]{\underset{^{\longrightarrow }}{LD }}  ^\mathrm{b} _{\Q, \mathrm{qc}}
(\overset{^\mathrm{g}}{} \smash{\widehat{\D}} _{\PP} ^{(\bullet)})$,
$Y$ un sous-schéma de $P$.
Le $F$-complexe $\E$ se dévisse sur $Y$ en $F$-isocristaux surconvergents si et seulement si, pour tout
sous-schéma $Y '\subset Y$ tel qu'il existe un diviseur $T'$ de $P$ satisfaisant
$Y ' = \smash{\overline{Y}} ' \setminus T '$, $\E$ se dévisse sur $Y '$ en $F$-isocristaux surconvergents.
\end{prop}

\begin{proof}
  La suffisance est tautologique. Traitons à présent la nécessité.
  Il existe
  un recouvrement fini ouvert $(Y _l) _l$ de $Y$,
des diviseurs $T _l$ de $P$ tels que, pour tout $l$,
$Y _l = \smash{\overline{Y}} _l \setminus T _l$ et $\E$ soit dévissable sur $Y _l$ en $F$-isocristaux surconvergents.
Soit maintenant $Y '\subset Y$ un sous-schéma tel qu'il existe un diviseur $T'$ de $P$ satisfaisant
$Y ' = \smash{\overline{Y}} ' \setminus T '$.
Par \ref{indYdev-lemm1}, $\E$ l'est alors sur $Y' \cap Y _l$
(comme $Y _l = Y \setminus T _l$, on remarque que $Y' \cap Y _l=\smash{\overline{Y}}' \setminus (T' \cup T _l)$).
On en déduit, via \ref{EdagtTT'}.\ref{EdagtTT'i}), \ref{EdagtTT'ii}), que $(\hdag T _l) (\E)$ l'est sur $Y'$.
Par \ref{lemm-devloc}, $(\hdag \cap _l T _l) (\E)$ l'est donc sur $Y'$. Comme $Y' \setminus (\cap _l T _l)=Y'$
(et donc $\R \underline{\Gamma} ^\dag _{Y'} (\E) \riso \R \underline{\Gamma} ^\dag _{Y'}( (\hdag \cap _l T _l) (\E) $),
cela équivaut à dire que $\E$ l'est sur $Y'$.
\hfill \hfill \qed \end{proof}

On déduit de \ref{devloc} les deux corollaires suivants.

\begin{coro}
\label{coro1devloc}
  Soient $\E \in F \text{-} \smash[b]{\underset{^{\longrightarrow }}{LD }}  ^\mathrm{b} _{\Q, \mathrm{qc}}
(\overset{^\mathrm{g}}{} \smash{\widehat{\D}} _{\PP} ^{(\bullet)})$,
$Y' \subset Y$ deux sous-schémas de $P$.
Si le $F$-complexe $\E$ se dévisse sur $Y$ en $F$-isocristaux surconvergents alors $\E$ l'est sur $Y'$.
\end{coro}

\begin{coro}
\label{coro2devloc}
  Soient $\E \in F \text{-} \smash[b]{\underset{^{\longrightarrow }}{LD }}  ^\mathrm{b} _{\Q, \mathrm{qc}}
(\overset{^\mathrm{g}}{} \smash{\widehat{\D}} _{\PP} ^{(\bullet)})$,
$Y', Y$ deux sous-schémas de $P$.
Le $F$-complexe $\E$ est dévissable sur $Y \cup Y'$ en $F$-isocristaux surconvergents si et seulement s'il l'est sur
$Y$ et $Y'$.
\end{coro}

\begin{prop}
  \label{gammaYdevup}
  Soient $\E \in F \text{-} \smash[b]{\underset{^{\longrightarrow }}{LD }}  ^\mathrm{b} _{\Q, \mathrm{qc}}
(\overset{^\mathrm{g}}{} \smash{\widehat{\D}} _{\PP} ^{(\bullet)})$,
$Y, Y_1, Y_2$ trois sous-schémas de $P$.

Si $\R \underline{\Gamma} ^\dag _{Y_1} (\E) \riso
\R \underline{\Gamma} ^\dag _{Y_2} (\E)$
alors $\E$ est dévissable sur $Y\cap Y _1$ en $F$-isocristaux surconvergents
si et seulement s'il l'est sur $Y \cap Y_2$.
\end{prop}

\begin{proof}
Supposons $\E$ dévissable sur $Y\cap Y _1$ en $F$-isocristaux surconvergents.
Avec \ref{gammayY'}, on obtient :
$\R \underline{\Gamma} ^\dag _{(Y \setminus Y_1) \cap Y _2} (\E) \riso
\R \underline{\Gamma} ^\dag _{(Y \setminus Y_1) \cap Y _1} (\E) = 0.$
Comme $ Y \cap Y _2 = (Y \cap Y _2 \cap Y _1 )\cup (Y  \cap Y _2 \setminus Y _1 )$,
grâce à \ref{coro1devloc} et \ref{coro2devloc}, on conclut alors que $\E$ est
dévissable sur $Y\cap Y _2$ en $F$-isocristaux surconvergents.
\hfill \hfill \qed \end{proof}

\begin{defi}
\label{def-dev-genU}
Soient $T$ un diviseur de $P$,
$\E \in F \text{-} \smash[b]{\underset{^{\longrightarrow }}{LD }}  ^\mathrm{b} _{\Q, \mathrm{qc}}
(\overset{^\mathrm{g}}{} \smash{\widehat{\D}} _{\PP} ^{(\bullet)} (T ))$
($\subset F \text{-} \smash[b]{\underset{^{\longrightarrow }}{LD }}  ^\mathrm{b} _{\Q, \mathrm{qc}}
(\overset{^\mathrm{g}}{} \smash{\widehat{\D}} _{\PP} ^{(\bullet)} )$.
Le $F$-complexe $\E$ {\og se dévisse en $F$-isocristaux surconvergents \fg}
si
$\E$ se dévisse sur $P$ (ou sur $P \setminus T$ d'après \ref{EdagtTT'}.\ref{EdagtTT'ii})) en $F$-isocristaux surconvergents.

La sous-catégorie pleine de
$F \text{-} \smash[b]{\underset{^{\longrightarrow }}{LD }}  ^\mathrm{b} _{\Q, \mathrm{qc}}
(\overset{^\mathrm{g}}{} \smash{\widehat{\D}} _{\PP} ^{(\bullet)} (T ))$ des
$F$-complexes dévissables
en $F$-isocristaux surconvergents
sera noté
$F \text{-} \smash[b]{\underset{^{\longrightarrow }}{LD }}  ^\mathrm{b} _{\Q, \textrm{dév}}
(\overset{^\mathrm{g}}{} \smash{\widehat{\D}} _{\PP} ^{(\bullet)} (T ))$.
\end{defi}

Dans la suite, $T$ désignera un diviseur de $P$.

\begin{prop}
\label{qcdevtriang}
La catégorie
  $F \text{-} \smash[b]{\underset{^{\longrightarrow }}{LD }}  ^\mathrm{b} _{\Q, \textrm{dév}}
(\overset{^\mathrm{g}}{} \smash{\widehat{\D}} _{\PP} ^{(\bullet)} (T ))$
est une sous-catégorie pleine triangulée de
$F \text{-} \smash[b]{\underset{^{\longrightarrow }}{LD }}  ^\mathrm{b} _{\Q, \textrm{qc}}
(\overset{^\mathrm{g}}{} \smash{\widehat{\D}} _{\PP} ^{(\bullet)} (T ))$.
\end{prop}
\begin{proof}
  Cela résulte de \ref{EE'dev}.
\hfill \hfill \qed \end{proof}

\begin{defi}
  \label{rema-def-dev}
Soient $n$ un entier naturel,
$\E \in F \text{-} \smash[b]{\underset{^{\longrightarrow }}{LD }}  ^\mathrm{b} _{\Q, \textrm{qc}}
(\overset{^\mathrm{g}}{} \smash{\widehat{\D}} _{\PP} ^{(\bullet)} (T ))$.
  On définit la notion de {\og $n$-dévissabilité en $F$-isocristaux surconvergents\fg} par récurrence sur $n$ de la
  manière suivante :
  \begin{itemize}
    \item Le $F$-complexe $\E$ est $0$-dévissable en $F$-isocristaux surconvergents
    s'il existe un sous-schéma lisse $Y$ de $P$, un diviseur $T'\supset T$ de $P$ tels que
  $Y=\overline{Y} \setminus T'$,
  $\R \underline{\Gamma} ^\dag _{Y}(\E) \riso \E$ et
  les espaces de cohomologie de $\E$ soient associés à des
  $F$-isocristaux surconvergents sur $Y$ ;
  \item Si $n \geq 1$, le $F$-complexe $\E$ est $n$-dévissable en $F$-isocristaux surconvergents s'il existe
  un triangle distingué de la forme :
  $$ \E _0 \rightarrow \E \rightarrow \E _{n-1} \rightarrow \E _0 [1],$$
  où $\E _0$ est $0$-dévissable en $F$-isocristaux surconvergents
  et $\E _{n-1}$ est $n-1$-dévissable en $F$-isocristaux surconvergents.
  \end{itemize}
Enfin, {\og $\E$ est $\infty$-dévissable en $F$-isocristaux surconvergents\fg}
  s'il existe un entier $n \geq 0$ tel que $\E$ soit $n$-dévissable en $F$-isocristaux surconvergents.
\end{defi}
La définition \ref{rema-def-dev} paraît asymétrique. La proposition \ref{n1n2dev} implique qu'il n'en est rien.
\begin{prop}
  \label{n1n2dev}
  Soient $n _1, n _2\geq 0$ deux entiers,
  $\E _{n _1} \rightarrow \E _{n _2} \rightarrow \E  \rightarrow \E _{n _1} [1]$ un triangle distingué de
$F \text{-} \smash[b]{\underset{^{\longrightarrow }}{LD }}  ^\mathrm{b} _{\Q, \textrm{qc}}
(\overset{^\mathrm{g}}{} \smash{\widehat{\D}} _{\PP} ^{(\bullet)} (T ))$
 tel que $\E _{n _1}$ soit $n _1$-dévissable en $F$-isocristaux surconvergents,
 $\E _{n _2}$ soit $n _2$-dévissable en $F$-isocristaux surconvergents. Alors $\E $ est $n _1 + n_2 + 1$-dévissable
 en $F$-isocristaux surconvergents.
\end{prop}

\begin{proof}
  On procède par récurrence sur $n _2$. Si $n _2 =0$, c'est tautologique. Supposons $n _2 \geq 1$.
  Par définition, on dispose d'un triangle distingué
  $\E _{n _2} \rightarrow \E _{n _2 -1} \rightarrow \E _0 \rightarrow \E _{n _2} [1]$, où
  $\E _{n _2 -1}$ est $n _2 -1$-dévissable en $F$-isocristaux surconvergents et
  $\E _0$ est $0$-dévissable en $F$-isocristaux surconvergents. On obtient par composition un morphisme
  $\E _{n _1} \rightarrow \E _{n _2 -1}$. En notant $\E _{n _1 +n _2}$ son cône, par hypothèse de récurrence,
  $\E _{n _1 +n _2}$ est $n _1 + n_2$-dévissable en $F$-isocristaux surconvergents.
  L'axiome de l'octaèdre nous donne un triangle distingué de la forme :
  $\E  \rightarrow \E _{n _1 +n _2} \rightarrow \E _0 \rightarrow \E  [1]$. D'où le résultat.
\hfill \hfill \qed \end{proof}

\begin{prop}
\label{ppetitetrgl}
Un $F$-complexe $\E\in F \text{-} \smash[b]{\underset{^{\longrightarrow }}{LD }}  ^\mathrm{b} _{\Q, \textrm{qc}}
(\overset{^\mathrm{g}}{} \smash{\widehat{\D}} _{\PP} ^{(\bullet)} (T ))$
est dévissable en $F$-isocristaux surconvergents si et seulement s'il
  est $\infty$-dévissable en $F$-isocristaux surconvergents.
De plus, $F \text{-} \smash[b]{\underset{^{\longrightarrow }}{LD }}  ^\mathrm{b} _{\Q, \textrm{dév}}
(\overset{^\mathrm{g}}{} \smash{\widehat{\D}} _{\PP} ^{(\bullet)} (T ))$
est la plus petite sous-catégorie pleine triangulée de
$F \text{-} \smash[b]{\underset{^{\longrightarrow }}{LD }}  ^\mathrm{b} _{\Q, \textrm{qc}}
(\overset{^\mathrm{g}}{} \smash{\widehat{\D}} _{\PP} ^{(\bullet)} (T ))$
contenant les $F$-complexes $0$-dévissables en $F$-isocristaux surconvergents.
\end{prop}

\begin{proof}
  Cela résulte de \ref{qcdevtriang}.
\hfill \hfill \qed \end{proof}

On dispose en outre de
la sous-catégorie pleine
$F\text{-}D ^\mathrm{b} _{\textrm{dév}}
(\smash{\D} ^\dag _{\PP} (\hdag T) _\Q)$ de
$F\text{-}D ^\mathrm{b} _{\textrm{coh}}
(\smash{\D} ^\dag _{\PP} (\hdag T) _\Q)$
des $F$-complexes dévissables en $F$-isocristaux surconvergents (voir \cite[8.1.1]{caro_devissge_surcoh}).
La proposition qui suit indique que ces deux notions de dévissabilité se rejoignent.

\begin{prop}
\label{dev=qwcdev+coh}
Le foncteur $\underset{\longrightarrow}{\lim}$
induit une équivalence de catégories :
\begin{equation}
  \label{dev=qwcdev+cohlim}
\underset{\longrightarrow}{\lim} \ : \
F \text{-} \smash[b]{\underset{^{\longrightarrow }}{LD }}  ^\mathrm{b} _{\Q, \textrm{dév}}
(\overset{^\mathrm{g}}{} \smash{\widehat{\D}} _{\PP} ^{(\bullet)} (T ))
\cap
F \text{-} \smash[b]{\underset{^{\longrightarrow }}{LD }}  ^\mathrm{b} _{\Q, \textrm{coh}}
(\overset{^\mathrm{g}}{} \smash{\widehat{\D}} _{\PP} ^{(\bullet)} (T ))
\cong
F\text{-}D ^\mathrm{b} _{\textrm{dév}}
(\smash{\D} ^\dag _{\PP} (\hdag T) _\Q).
\end{equation}
\end{prop}
\begin{proof}
Soit $\E ^{(\bullet)} \in
F \text{-} \smash[b]{\underset{^{\longrightarrow }}{LD }}  ^\mathrm{b} _{\Q, \textrm{dév}}
(\overset{^\mathrm{g}}{} \smash{\widehat{\D}} _{\PP} ^{(\bullet)} (T ))
\cap
F \text{-} \smash[b]{\underset{^{\longrightarrow }}{LD }}  ^\mathrm{b} _{\Q, \textrm{coh}}
(\overset{^\mathrm{g}}{} \smash{\widehat{\D}} _{\PP} ^{(\bullet)} (T ))$.
Via l'équivalence
$\underset{\longrightarrow}{\lim}$ :
$F \text{-} \smash[b]{\underset{^{\longrightarrow }}{LD }}  ^\mathrm{b} _{\Q, \textrm{coh}}
(\overset{^\mathrm{g}}{} \smash{\widehat{\D}} _{\PP} ^{(\bullet)} (T ))
\cong
F\text{-}D ^\mathrm{b} _{\textrm{coh}}
(\smash{\D} ^\dag _{\PP} (\hdag T) _\Q)$,
on note $\E := \underset{\longrightarrow}{\lim}\,( \E ^{(\bullet)})$,
$X$ son support et $Y := X \setminus T$.
Par \ref{indYdev-lemm1},
$\E ^{(\bullet)}$ se dévisse en $F$-isocristaux surconvergents au-dessus d'une stratification
de la forme
$Y=\sqcup _{i= 0,\dots , r-1} Y _i$.
Avec \ref{remacaradev},
on peut dévisser $\E ^{(\bullet)}$ au-dessus de $Y _i$ en utilisant une stratification dont les strates
sont affines, lisses et munis de modèles idéaux.
En utilisant \ref{lemm-devdev=dev}, on obtient ainsi un dévissage de $\E$ au sens de
\cite[8.1.1]{caro_devissge_surcoh} (car les strates sont munies de modèles idéaux).
Le foncteur $\underset{\longrightarrow}{\lim}$ de \ref{dev=qwcdev+cohlim} est donc bien défini.
Comme sa pleine fidélité est immédiate, vérifions à présent son essentielle surjectivité.
Soient
$\E \in F\text{-}D ^\mathrm{b} _{\textrm{dév}}
(\smash{\D} ^\dag _{\PP} (\hdag T) _\Q)$,
$X$ son support et $Y := X \setminus T$.
Il existe $\E ^{(\bullet)}
\in F \text{-} \smash[b]{\underset{^{\longrightarrow }}{LD }}  ^\mathrm{b} _{\Q, \textrm{coh}}
(\overset{^\mathrm{g}}{} \smash{\widehat{\D}} _{\PP} ^{(\bullet)} (T ))$
tel que $\underset{\longrightarrow}{\lim} ( \E ^{(\bullet)})\riso \E$.
Dans la définition de la dévissabilité de \cite[8.1.1]{caro_devissge_surcoh},
on peut, grâce à \ref{lemm-devdev=dev} et \ref{remacaradev}, éviter l'hypothèse {\og $X _r$
(notation de \cite[8.1.1]{caro_devissge_surcoh})
est lisse\fg}.
Ainsi $\E ^{(\bullet)}$ se dévisse sur $Y$ en $F$-isocristaux surconvergents (sens de \ref{coro-caradev2}).
Comme $\E $ est à support dans $X$, par \ref{lemgammaYdevup},
$\E ^{(\bullet)}$ est donc dévissable en $F$-isocristaux surconvergents.
\hfill \hfill \qed \end{proof}

Dans \cite[3.17]{caro_surholonome}, pour toute $k$-variété $Y $,
nous avons défini la notion de
$F\text{-}\smash{\D} _Y$-module arithmétique surholonome
(resp. $F\text{-}\smash{\D} _Y$-complexe arithmétique surholonome).
En gros, un $F$-complexe cohérent et bornée de $\D$-modules arithmétiques
est surholonome si sa cohérence est stable par foncteurs duaux,
par images inverses extraordinaires par un morphisme lisse, par
foncteurs restrictions en dehors d'un diviseur
et par composition de tels foncteurs.  \'Enonçons à présent les deux conjectures ci-dessous.

\begin{conj}
\label{conj}
\begin{enumerate}
\item \label{conj1}
  Soit $Y $ une $k$-variété lisse.
  Pour tout $F$-isocristal $E$ surconvergent sur $Y$,
  $\sp _{Y+} (E)$ est un $F\text{-}\smash{\D} _Y$-module arithmétique surholonome.

  \item \label{conj2} On dispose de l'égalité :
$$F \text{-} \smash[b]{\underset{^{\longrightarrow }}{LD }}  ^\mathrm{b} _{\Q, \textrm{dév}}
(\overset{^\mathrm{g}}{} \smash{\widehat{\D}} _{\PP} ^{(\bullet)} (T ))
=
F \text{-} \smash[b]{\underset{^{\longrightarrow }}{LD }}  ^\mathrm{b} _{\Q, \textrm{dév}}
(\overset{^\mathrm{g}}{} \smash{\widehat{\D}} _{\PP} ^{(\bullet)} (T ))
\cap
F \text{-} \smash[b]{\underset{^{\longrightarrow }}{LD }}  ^\mathrm{b} _{\Q, \textrm{coh}}
(\overset{^\mathrm{g}}{} \smash{\widehat{\D}} _{\PP} ^{(\bullet)} (T )).$$
\end{enumerate}
\end{conj}

\begin{rema}
\begin{enumerate}
\item
En notant $U := P \setminus T$,
la conjecture \ref{conj}.\ref{conj1}) implique aussitôt \ref{conj}.\ref{conj2})
ainsi que les équivalences de catégories (notation de \cite[3.13]{caro_surholonome}):
$$F \text{-} \smash[b]{\underset{^{\longrightarrow }}{LD }}  ^\mathrm{b} _{\Q, \textrm{dév}}
(\overset{^\mathrm{g}}{} \smash{\widehat{\D}} _{\PP} ^{(\bullet)} (T ))
\cong
F\text{-}D ^\mathrm{b} _{\textrm{dév}}
(\smash{\D} ^\dag _{\PP} (\hdag T) _\Q)
\cong F\text{-}D ^\mathrm{b} _\mathrm{surhol} (\D _{U }).$$

\item  La conjecture \cite[8.1.6]{caro_devissge_surcoh} est équivalente à \ref{conj}.\ref{conj1}).
En effet,
avec les notations et hypothèses de \ref{conj}.\ref{conj1}),
on peut toujours dévisser $\sp _{Y+} (E)$ au-dessus d'une stratification dont les strates sont affines, lisses
et munis de modèles idéaux (voir \ref{remacaradev}).
Cela entraîne que la conjecture \cite[8.1.6]{caro_devissge_surcoh} implique \ref{conj}.\ref{conj1}).
La réciproque est une tautologie.

\item

Toutes les constructions (e.g. fonctions $L$ généralisant celles de \'Etesse et Le Stum de \cite{E-LS} ou \cite{E-LS2},
notion de poids, stabilité du poids par image directe etc.) et résultats de \cite[8]{caro_devissge_surcoh}
sont encore valables en remplaçant
$F\text{-}D ^\mathrm{b} _{\textrm{dév}}
(\smash{\D} ^\dag _{\PP} (\hdag T) _\Q)$
par
$F \text{-} \smash[b]{\underset{^{\longrightarrow }}{LD }}  ^\mathrm{b} _{\Q, \textrm{dév}}
(\overset{^\mathrm{g}}{} \smash{\widehat{\D}} _{\PP} ^{(\bullet)} (T ))$.
\end{enumerate}

\end{rema}

La proposition qui suit valide la conjecture \ref{conj}.\ref{conj2}).

\begin{prop}\label{courbeholodev}
  Si $P$ est de dimension au plus $1$ alors
  $$F \text{-} \smash[b]{\underset{^{\longrightarrow }}{LD }}  ^\mathrm{b} _{\Q, \textrm{dév}}
(\overset{^\mathrm{g}}{} \smash{\widehat{\D}} _{\PP} ^{(\bullet)} (T ))
=
F \text{-} \smash[b]{\underset{^{\longrightarrow }}{LD }}  ^\mathrm{b} _{\Q, \textrm{dév}}
(\overset{^\mathrm{g}}{} \smash{\widehat{\D}} _{\PP} ^{(\bullet)} (T ))
\cap
F \text{-} \smash[b]{\underset{^{\longrightarrow }}{LD }}  ^\mathrm{b} _{\Q, \textrm{coh}}
(\overset{^\mathrm{g}}{} \smash{\widehat{\D}} _{\PP} ^{(\bullet)} (T )).$$
\end{prop}

\begin{proof}
Cela résulte de l'holonomie des $F$-isocristaux surconvergents sur les courbes lisses
(voir \cite[4.3.4]{caro_courbe-nouveau}).
\hfill \hfill \qed \end{proof}

\bibliographystyle{smfalpha}
\bibliography{bib1}

\end{document}